%% file: defselmodesaim3.tex
\documentclass{ps}

\usepackage[latin1]{inputenc}
  \usepackage[T1]{fontenc}
\usepackage[greek,francais,english]{babel}

\usepackage[dvipsone]{color}
\usepackage[dvips]{graphicx}
\usepackage{amsmath}

\usepackage{latexsym}

\usepackage{amsfonts,amssymb}

\newcommand{\suptwo}[2]{\sup_{\substack{#1 \\ #2}}} 

\newcommand{\norm}[1]{\left\lVert #1 \right\rVert} 
\newcommand{\ps}[2]{\langle #1,#2 \rangle} 
\newcommand{\esp}[1]{\mathbb{E}\left[ #1 \right]} 
\newcommand{\prob}[2]{\mathbb{P}_{#1}\left[ #2 \right]} 
\newcommand{\dd}{\mathrm{d}} 
\newcommand{\rv}[1]{{\bf \{~#1~\}}} 

\DeclareMathOperator{\pen}{pen}
\DeclareMathOperator{\Vect}{Vect}

\begin{document}

\title{Model selection for quantum homodyne tomography}
\author{Jonas Kahn}
\address{Universit\'e
Paris-Sud 11, D\'epartement de Math\'ematiques
 B\^at 425, 91405 Orsay Cedex,
France; e-mail: jonas.kahn@math.u-psud.fr}
\subjclass{62G05, 81V80, 62P35}
\keywords{density matrix, model selection, pattern functions estimator,
penalized maximum likelihood estimator, penalized projection estimators,
quantum calibration, quantum tomography, wavelet estimator, Wigner function.}

\setlength{\parindent}{0pt}
\setlength{\parskip}{1mm}

\begin{resume}
Nous nous intéressons à un problème de statistique non-paramétrique issu de la physique, et plus précisément à la tomographie quantique, c'est-à-dire la détermination de l'état quantique d'un mode de la lumière via une mesure                           homodyne. Nous appliquons plusieurs procédures de sélection de modèles: des estimateurs par projection pénalisés, où on peut utiliser soit des fonctions motif, soit des ondelettes, et l'estimateur du maximum de vraisemblance pénalisé. Dans chaque cas, nous obtenons une inégalité oracle. Nous prouvons également une vitesse de convergence polynomiale pour ce problème non-paramétrique, pour les estimateurs par projection. Nous appliquons ensuite des idées à la calibration d'un photocompteur, l'appareil dénombrant le nombre de photons dans un rayon lumineux. Le problème mathématique se réduit dans ce cas à un problème non-paramétrique à données manquantes. Nous obtenons à nouveau des inégalités oracle, qui nous assurent des vitesses de convergence d'autant meilleures que le photocompteur est bon.
\end{resume}

\begin{abstract}
 This paper deals with a non-parametric problem coming from physics, namely quantum tomography. That consists in determining the quantum
 state of a mode of light through a homodyne measurement.  We 
apply several model selection procedures: penalized projection estimators,
where we may use pattern functions or wavelets, and penalized maximum
likelihood estimators. In all these cases, we get oracle inequalities. In the
former we also have a polynomial  rate of convergence for the non-parametric
problem. We finish the paper with applications of similar ideas to the
calibration of a photocounter, a measurement apparatus counting the number of
photons in a beam. Here the mathematical problem reduces similarly to a
non-parametric missing
data problem. We again get oracle inequalities, and better speed if the
photocounter is good.
\end{abstract}

\maketitle

\setlength{\parindent}{0pt}
\setlength{\parskip}{2mm}

\section{Introduction}
  
  Quantum mechanics introduces intrinsic randomness in physics: the result of a
measurement, or any macroscopic interaction, on a physical system is not deterministic. Therefore, a host of
statistical problems can stem from it. Some are (almost) specifically quantum,
notably any question about which measurement yields the maximum information, or
whether simultaneously measuring $n$ samples is more efficient than measuring
them sequentially \cite{gill-2001-36}. However, once we have chosen the
measurement we
carry out on our physical system, we are left with an entirely classical
statistical problem. This paper aims at applying  model selection
methods \emph{\`a la} Birg\'e-Massart to one such instance, which is of
interest both practical, as physicists use this measurement quite often (the
underlying physical system is elementary; it is the particle with one degree of
freedom), and mathematical, as it yields a nonparametric inverse problem with uncommon features.

Moreover, as this classical problem stemming from quantum mechanics could be
seen as an easy introduction to the subject to classical statisticians, we have
added more general notions on quantum statistics at the beginning of the appendix. The interested reader can get further acquaintance with these
concepts through the textbooks  \cite{Helstrom} and \cite{Holevo} or the review
article \cite{Barndorff-Nielsen&Gill&Jupp}.

   More precisely, the problem we are interested in is quantum homodyne tomography. As an
aside, we apply the results we get to the calibration of a photocounter,
using a quantum tomographer as a tool.  The word ``Homodyne'' refers to the
experimental technique used for this measurement, first implemented in
\cite{Smithey}, where the state of one mode of electromagnetic
radiation, that is a pulse of laser light at a given frequency, is probed using
a reference laser beam at the same (``homo'') frequency. And ``Tomography'' is used because one of the physicists' favourite representations of the state, the Wigner function, can be recovered from the data by inverting a Radon transform. 

  Mathematically, our data are samples from a probability distribution
$p_\rho$  on $ \mathbb{R} \times [0,\pi] $. From this data, we want to recover
the ``density operator'' $\rho$ of the system. This is the most common representation
of the state, that is a mathematical object which encodes all the information
about the system. Perfect knowledge of the state means knowing how the system
will evolve and the probability distribution of the result of any measurement
we might carry out on the system. And these laws of evolution and
measurement can be expressed naturally enough within the density operator
framework (see Appendix). The density operator is a 
non-negative trace-one self-adjoint
operator $\rho$ on $ L^2( \mathbb{R} )$ (in our particular case). We know the
linear transform $ \mathcal{T} $ which takes $\rho$ to $p_\rho$ and can make it
explicit in particular bases such as the Fock basis. We may also settle for the
Wigner function $W$, another representation of the state. That is a
two-dimensional real function with integral one, and $p_\rho$ is the Radon transform of $W$. 

  The first reconstruction methods used the Wigner function as an intermediate representation:
after collecting the data in histograms and smoothing, one inverted the Radon
transform to get an estimate of $W$. This smoothing, however, introduces hard-to-control bias. Using the
pattern functions (bidual bases, in fact) introduced in
\cite{D'Ariano.0}, consistency of linear estimators of the density operator was
proved in  \cite{Gu.Gill.}. There were also similar
results for sieved maximum likelihood estimators. Then, 
a sharp adaptive estimator for the Wigner function was devised in \cite{But.Gu.Art.}, and this even if there is noise in the measurement (see subsection \ref{noise}). 

  In this paper, we devise penalized estimators that fulfill oracle-type
inequalities among the $L^2$-projections on submodels, analyze the penalized
maximum likelihood estimator and apply these estimators to the calibration of a
photocounter. Notice that all these results are derived for finite samples
(all the previous works considered only the asymptotic regime). We have mainly worked under the idealized hypothesis where there is no noise, however.

  The appendix is not logically necessary for the article. We
have            inserted it for background and as an invitation to this field.
It first features a general introduction to quantum
statistics with a public of classical statisticians in mind. We then describe what 
quantum homodyne tomography precisely is.
This latter subsection is largely
based on \cite{But.Gu.Art.}.

 Section \ref{Problème} formalizes the statistical problem at hand, with no
need of the appendix, except the equations explicitly referred to
therein.

 Section \ref{projection} aims at devising a model selection procedure to
choose between $L^2$-projection estimators. We first give general theorems
(\ref{deter} and \ref{alea}) leading to oracle-type inequalities for
hard-thresholding estimators. We then apply them to two bases. One is the Fock
basis and the corresponding pattern functions physicists have used for a while.
For it we also prove a polynomial convergence rate for any state with finite
energy. The other is a wavelet basis for the Wigner function. We finish with a
short subsection describing what changes are entailed by the presence of noise. Especially, we do not need to adapt our theorems if the noise is low enough, as long as we change the dual basis.

 Section \ref{EMV} similarly applies a classical theorem (\ref{Massart}) to solve the question of which (size of) model is best to use a maximum likelihood estimator on. 

  Section \ref{photocounter} switches to the determination of a kind of
measurement apparatus (and not any more on the state that is sent in) using a
known state and this same tomographer that was studied in the previous
sections. The law of our samples are then very similar and we apply the same
type of techniques (penalized projection and maximum likelihood estimators).
The fact that the POVM (mathematical modelling of a measurement) is a projective measurement (see Appendix) enables us to work with $L^1$-operator norm, however.

\section{The mathematical problem}
\label{Problème}

We now describe the mathematical problem at hand.

We are given $n$ independent identically distributed random variables  $Y_i =
(X_i,\phi_i)$ with   density $p_\rho$ on $[0,\pi)\times \mathbb{R}$.

This data is the result of a measurement on a physical system. Now the ``state'' of a system is described by a mathematical object, and there are two favourites for physical reasons: one is the \emph{density operator} $\rho$, the other is the \emph{Wigner function} $W_{\rho}$. We describe them below.

Therefore we are not actually interested in $p_\rho$, but rather in $W_{\rho}$ or (maybe preferably) $\rho$. The probability distribution $p_{\rho}$ of our samples can be retrieved if we know either $\rho$ or $W_{\rho}$.

In other words we aim at estimating as precisely as possible $\rho$ or $W_{\rho}$ from the data $\{Y_i\}$. By `` as precisely as possible'', we mean that with a suitable notion of distance, we shall minimize $\esp{d(\rho,\hat{\rho})}$. Our choice of distance will be partly dictated by mathematical tractability. 

\smallskip 

We now briefly explain what $W_{\rho}$ and $\rho$ stand for.

\begin{figure}[ht]
\centering
\includegraphics[height=0.2\textheight]{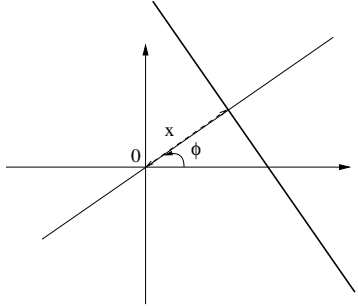}
\caption{\label{Radon} The value of $p_\rho$ at $(x,\phi)$ is the integral of the Wigner function over the bold line}
\end{figure}

The Wigner function $W_{\rho}: \mathbb{R}^2 \rightarrow \mathbb{R}$ is the inverse Radon transform of $p_{\rho}$. In fact we would rather say that $p_{\rho}$ is the Radon transform of $W_{\rho}$. 
 Explicitly:
\begin{eqnarray*}
p_{\rho}(x,\phi) & = & \int_{-\infty}^{\infty} W(x\cos\phi + y \sin \phi,x \sin
\phi - y\cos \phi)dy.
\end{eqnarray*} 
Figure \ref{Radon} might be of some help.
An important remark is that  the Wigner function is not a probability density, but only a quasi-probability density: a function with integral $1$, but that may be negative at
places. However its Radon transform is a true probability density, as it is $p_{\rho}$.

Retrieving $W_{\rho}$ from $P_{\rho}$ then amounts to inverting the Radon transform, hence the name of tomography: that is the same mathematical problem as with the brain imagery technique called Positron Emission Tomography.  

As for $\rho$, this is a density operator on the Hilbert space $L^2(\mathbb R)$, that is a \emph{self-adjoint positive} operator with \emph{trace $1$}. We denote the set of such operators by $ \mathcal{S}(L^2(\mathbb{R}))$. 
There is a linear transform $\mathbf T$ that takes $\rho$ to $p_{\rho}$. We give it explicitly using a basis of $ L^2(\mathbb{R})$ known as the \emph{Fock basis} This orthonormal basis, which has many nice physical properties, is defined by: 
\begin{eqnarray}
\label{hermite}
  \psi_k(x)&=&H_k(x)e^{-x^2/\!2}
\end{eqnarray}  
  where $H_k$ is the $k$th Hermite polynomial normalized such that $\norm{\psi_k}_2=1$. The matrix entries of $\rho$ in this basis are $\rho_{j,k} = \langle \psi_j , \rho \psi_k \rangle$. Then $\mathbf T$ can be written:
\begin{eqnarray*}
\bold{T}: \mathcal{S}(L^2(\mathbb{R})) & \longrightarrow     & L^1(\mathbb{R}\times [0,\pi]) \\
            \rho   & \mapsto  &   \left(p_{\rho}:(x,\phi)\mapsto
\sum_{j,k=0}^{\infty}\rho_{j,k}\psi_j(x)\psi_k(x)e^{-i(j-k)\phi}\right).
\end{eqnarray*}

Notice that as we have defined precisely the set of possible $\rho$, this mapping yields the set of possible $p_{\rho}$ and $W_{\rho}$.

The relations between $\rho$, $W_{\rho}$ and $p_{\rho}$ are further detailed in subsection \ref{QHD}.

\smallskip

Anyhow we may now state our problem as consisting in inverting either the Radon transform or $\mathbf T$ from empirical data.

  This is a classical problem of non-parametric statistics, that we want to treat non-asymptotically. We then take estimators based on a \emph{model}, that is a subset of the operators on $L^2(\mathbb{R})$, or equivalently of the two-dimensional real functions. These models are usually  vector spaces, which may not be the domain of the object to be estimated. To choose a candidate within a given model, there are different methods, two of which we study, projection estimators and maximum likelihood estimators. Once we have a candidate within each model, we then use model selection methods to choose (almost) the best.

  We  first study projection estimators, for which the most convenient distance comes from the $L^2$ norm
\begin{align*}
\|\tau\|_2 & = \sqrt{\sum | \lambda_i(\tau)|^2}  = \sqrt{\sum_{j,k} |\tau_{j,k}|^2},
\end{align*} 
where the $ \lambda_i$ are the eigenvalues of $ \tau$, and
the second equality holds for $\tau$ written in any orthonormal basis. Notice that there is an isometry (up to a constant) between the space of density operators with $L^2$-operator norm and the space of Wigner functions with $L^2$-Lebesgue norm, that is:
\begin{equation*}
\| W_\rho -W_\tau\|_2^2 = \int\!\!\int \left|W_\rho(q,p) - W_\tau(q,p)
\right|^2 dp\, dq = \frac 1{2\pi} \|\rho - \tau\|^2_2. 
\end{equation*}

  For maximum likelihood estimators, we have to make do with the weaker
Hellinger distance (see later \eqref{defhell}) on $L^{1}\big(\mathbb{R} \times [0,\pi]\big)$, to which $p_{\rho}$ belongs.

\section{Projection estimators}
\label{projection}

  In this section, which owes much to \cite{Mas}, we apply penalization
procedures to projection estimators. The first subsection explains that we want
to obtain oracle-type inequalities. In the second we obtain a general
inequality where the left-hand side corresponds to an oracle inequality, and
where the remainder term in the right-hand side depends on the penalty and on
the large deviations of empirical coefficients. The two following subsections give two ways to choose the penalty term large enough for this remainder term to be small enough. In section \ref{deterpen} this penalty is deterministic. We design it and prove that it is a ``good choice'' by keeping Hoeffding's inequality in mind. In section \ref{ranpen}, the penalty is random, and designed by taking Bernstein's inequality into account. 

We next express these theorems in terms of two specific bases. For the Fock basis, we obtain
polynomial worst-case convergence rates, using the structure of states. For a
wavelet basis, we notice we obtain a usual estimator in classical tomography.
We finish by saying what can be done if there is noise, that is
(mainly) convolution of
the law of the sample by a gaussian.  We multiply the Fourier transform of the dual basis with the inverse of the Fourier transform of the gaussian, and as long as we still have well-defined functions, and we can re-use our theorems without changes.

\subsection{Aim of model selection}
\label{31}
  
  Let's assume we are given a (countable) $L^2$-basis $(e_i)_{i\in \mathcal{I}}$ of a
space in which $ \mathcal{S}( L^2(\mathbb{R}) )$ is included  (typically $ \mathcal{T}( L^2(\mathbb{R}) )$, the trace-class operators on $ L^2(\mathbb{R}) $). We may then try and find the coefficients of $\rho$ in this basis. The natural way to do so is to find a dual basis $(f_i)_{i\in \mathcal{I}}$ such that $\langle \textbf{T}(e_i) , f_j \rangle = \delta_{i,j}$ for all $i$ and $j$. Then, if $\rho= \sum_i \rho_i e_i$ we get $\langle p_\rho, f_i \rangle = \rho_i$ for all $i$. And if the $f_i$ are well enough behaved, then $\frac1{n}\sum_{k=1}^n f_i(X_k,\phi_k)=\hat{\rho_i}$  tends to $\rho_i$ by the law of large numbers.
  
  Now if we took $\sum_i \hat{\rho}_i e_i$ as an estimator of $\rho$, we would
have an infinite risk as the variance would be infinite. We must therefore
restrict ourselves to models $m \in \mathcal{M}$, that is $\Vect\left(e_i, i\in m
\right)$, where $m$ is a finite set, and $\mathcal{M}$ is a set of models (we might
take $\mathcal{M}$ smaller than the set of all finite sets of $\mathbb{N}$).  
  
  We may then write the loss as 
\begin{equation*}
\norm{\hat{\rho}_m-\rho}^2=  \sum_{i \not\in m} |\rho_i|^2 + \sum_{i \in
m} |\rho_i - \hat{\rho}_i|^2
\end{equation*}  
  where the first term is a bias (modelling error) and the second term is an estimation error. The risk would have this expression:  
\begin{equation*}
  \esp{\norm{\hat{\rho}_m-\rho}^2} = \sum_{i \not\in m} |\rho_i|^2 + \sum_{i
\in m} \esp{|\rho_i - \hat{\rho}_i|^2}
\end{equation*}
  where the expectation is taken with respect to $p_{\rho}$, since $\hat{\rho_i}$ depends on the $(X_k,\phi_k)$.
  
  If we use an arbitrary model $m$, we  probably have not
have struck a good balance between the bias term and the variance term. The
whole point of penalisation is to have a data-driven procedure to choose the
``best'' model. We are aiming at choosing the model with (almost) the lowest error. We would dream of obtaining:
\begin{equation*}
  \hat{m} = \arg\inf_{m \in \mathcal{M}} \norm{\hat{\rho}_m-\rho}^2.
\end{equation*}
   
   That is of course too ambitious. Instead, we shall obtain the following kind of bound, called an oracle inequality:
\begin{eqnarray}
\label{oracle}
  \esp{\left\{\norm{\hat{\rho}_{\hat{m}}-\rho}^2 -\left( C \inf_{m \in \mathcal{M}}\left(
d^2(\rho, m) + \pen(m) \right) \right)\right\} \vee 0} & \leq & \epsilon_n
\end{eqnarray}   
   where $d^2(\rho,m)$ is the bias of the model $m$, C is some constant,
independent of $\rho$, $\pen(m)$ is a \emph{penalty} associated to the model
$m$ (the bigger the model, the bigger the penalty) and $\epsilon_n$ depends only on
$n$ the number of observations, and goes to $0$ when $n$ is going to
infinity. We shall try to take the penalty of the order of the variance of the model.   

Notice that we have given in \eqref{oracle} an unusual form of oracle inequality. These inequalities are more often written as
\begin{eqnarray*}
  \esp{\norm{\hat{\rho}_{\hat{m}}-\rho}^2} & \leq & \left( C \inf_{m \in \mathcal{M}}\left(
d^2(\rho, m) +\esp{ \pen(m)} \right) \right) + \epsilon_n.
\end{eqnarray*}
Our form implies the latter.

  The strategy is the following: 
  
  First, rewrite the projection estimators as \emph{minimum contrast
estimators}, that is minimizers of a function (called the \emph{empirical contrast}
function, and written $\gamma_n$), which is the same for all models. We 
also demand that, for any $m$, this empirical contrast function converges to a \emph{contrast} function $\gamma$, the minimizer in $m$ of which is the projection of $\rho$ on $m$. 
  
  Second, find a penalty function that overestimates with high enough
probability $(\gamma - \gamma_n)(\hat{\rho}_m)$ for all $m$ simultaneously. Use of concentration inequalities is pivotal at this point.
  
 The   next section makes all this more explicit.
  
\subsection{Risk bounds and choice of the penalty function} 

  First we notice that the minimum of 
\begin{equation*}
\begin{split}
  \gamma(\tau) & =  \norm{\tau}^2 - 2 \langle \tau, \rho \rangle \\
           & =  \norm{\rho - \tau}^2 - \norm{\rho}^2
\end{split}
\end{equation*}
over a model $m$ is attained at the projection of $\rho$ on $m$. Moreover 
\begin{equation*}
  \gamma_n(\tau)= \norm{\tau}^2 - 2 \sum_{i} \frac1{n} \sum_{k=1}^n \tau_i f_i(X_k,\phi_k)
\end{equation*}
converges in probability to $\gamma$ for any $m$ (and all $\tau$ such that $\norm{\tau}=1$
simultaneously), as there is only a finite set of $i$ such that $\tau_i \neq 0$ for $\tau \in m$.

Now the minimum of $\gamma_n$ over $m$ is attained by 
\begin{equation*}
  \tau = \sum_{i \in m} \frac1{n} \sum_{k=1}^n f_i(X_k,\phi_k) e_i.
\end{equation*}

So we have succeeded in writing projection estimators as minimum contrast estimators. We then define our final estimator by:
\begin{equation*}
  \hat{\rho}^{(n)} = \hat{\rho}_{\hat{m}}
\end{equation*}
with
\begin{equation*}
\hat{m} = \arg\min_{m \in \mathcal{M}} \gamma_n(\hat{\rho}_m) + \pen_n(m)
\end{equation*}
where $\pen_n$ is a suitably chosen function depending on n, $m$ and possibly the data.

We then get, for any $m$, for any $\tau_m \in m$,
\begin{eqnarray}
\label{suite}
  \gamma_n(\hat{\rho}^{(n)}) + \pen_n(\hat{m}) \leq \gamma_n(\hat{\rho}_{m}) + \pen_n(m) \leq  \gamma_n(\tau_m) +
\pen_n(m).
\end{eqnarray}
What's more, for any $m$, for any $\tau_m \in m$,
\begin{eqnarray}
\label{gan}
\gamma_n(\tau_m) & = & \norm{\rho-\tau_m}^2 - \norm{\rho}^2 - 2 \nu_n(\tau_m)
\end{eqnarray}
with 
\begin{eqnarray*}
\nu_n(\tau)& =  & \langle \tau, \rho \rangle - \sum_i \sum_{k=1}^n \tau_i f_i(X_k,\phi_k) \\
          & = &  \sum_{i \in m} \tau_i (\rho_i - \hat{\rho}_i) + \sum_{i \not\in
m} \tau_i \rho_i.
\end{eqnarray*}

Putting together (\ref{suite}) and (\ref{gan}), we get, for all $m$ and $\tau_m \in m$:
\begin{eqnarray*} 
  \norm{\hat{\rho}^{(n)} - \rho}^2 \leq \norm{\tau_m - \rho}^2 + 2 \nu_n(\hat{\rho}^{(n)} - \tau_m) +
\pen_n(m) - \pen_n(\hat{m}).
\end{eqnarray*}

We then want to take penalties big enough to dominate the fluctuations $\nu_n$. Some manipulations will make this expression more tractable. 
First we bound $\nu_n(\hat{\rho}^{(n)} - \tau_m)$ by $\norm{\hat{\rho}^{(n)} - \tau_m}\chi_n(m \cup \hat{m})$, with
\begin{eqnarray*}
  \chi_n(m) & = & \suptwo{\tau \in m}{\norm{\tau}=1}\nu_n(\tau).
\end{eqnarray*}

Now the triangle inequality gives $\norm{\hat{\rho}^{(n)} - \tau_m} \leq \norm{\hat{\rho}^{(n)} - \rho} + \norm{\rho - \tau_m}$, so that:
\begin{eqnarray*}
\norm{\hat{\rho}^{(n)} - \rho}^2 & \leq & \norm{\rho- \tau_m}^2 + 2 \chi_n(m \cup
\hat{m})\norm{\rho-\hat{\rho}^{(n)}}+2 \chi_n(m \cup\hat{m})\norm{\rho-\tau_m} - \pen_n(\hat{m}) +
\pen_n(m). 
\end{eqnarray*}

For all $\alpha>0$, the following holds: 
\begin{eqnarray}
\label{abaabb} 
2 ab & \leq  & \alpha a^2 + \alpha^{-1} b^2
\end{eqnarray}

Using this twice, we get, for all $\epsilon>0$:
\begin{eqnarray*}
\frac{\epsilon}{2+\epsilon}\norm{\rho-\hat{\rho}^{(n)}}^2 & \leq & \left( 1 + \frac{2}{\epsilon}\right)
\norm{\rho - \tau_m}^2 + (1+\epsilon) \chi_n^2(m\cup \hat{m}) -\pen_n(\hat{m}) +\pen_n(m).
\end{eqnarray*}

Noticing that $\chi_n(m\cup \hat{m}) \leq \chi_n(m) + \chi_n(\hat{m})$ and putting our estimate of the error in the left-hand side:

\begin{eqnarray*}
\frac{\epsilon}{2+\epsilon}\norm{\rho-\hat{\rho}^{(n)}}^2 - \left\{\left( 1 + \frac{2}{\epsilon}\right)
\norm{\rho - \tau_m}^2 + 2 \pen(m)\right\} & \leq & (1+\epsilon) (\chi_n^2(\hat{m}) +
\chi_n^2(m)) -\pen_n(\hat{m}) - \pen_n(m).
\end{eqnarray*}
Now what we want to avoid is that our penalty is less than the fluctuations, so we separate this event and take its expectation:
\begin{multline}
\label{esperancereelle2}
 \esp{\left\{\frac{\epsilon}{2+\epsilon}\norm{\rho-\hat{\rho}^{(n)}}^2 - \left( \left( 1 +
\frac{2}{\epsilon}\right) \norm{\rho - \tau_m}^2 + 2 \pen_n(m) \right)\right\}\vee 0}
\\ 
\begin{aligned} &  \leq
\esp{\left\{(1+\epsilon)(\chi_n^2(\hat{m})+\chi_n^2(m))-\pen(\hat{m})-\pen(m)\right\}\vee 0} 
 \\
 &
 \leq  2 \esp{\sup_m \left\{(1+\epsilon)\chi_n^2(m) - \pen(m)\right\}\vee 0}. 
\end{aligned}
\end{multline}

Thus stated, our problem is to take a   penalty large enough to make the right-hand
side negligible, that is vanishing like $1/n$.

We shall use this form of $\chi_n(m)$:
\begin{eqnarray*}
  \chi_n(m)  ~=~  \suptwo{(\tau_i)_{i\in m}}{\sum \tau_i^2=1} \sum_{i \in m} \tau_i
(\rho_i - \hat{\rho}_i)  ~=~ \sqrt{\sum_{i\in m} \left|\rho_i-\hat{\rho}_i
\right|^2}
\end{eqnarray*}
so that 
\begin{eqnarray}
\label{chi2}
\chi_n(m)^2  ~=~ \sum_{i\in m} \left|\rho_i-\hat{\rho}_i        \right|^2
          ~=~  \sum_{i\in m} \left|\rho_i-\frac1n\sum_{k=1}^n
f_i(x_k,\phi_k)       \right|^2.
\end{eqnarray} 

\subsection{Deterministic penalty}
\label{deterpen}

First we may try to craft a deterministic penalty. 
  
We plan to use Hoeffding's inequality, recalling that $\hat{\rho}_i$ is a sum
of independent variables:
\begin{lmm}{\bf: Hoeffding's inequality}
  Let $X_1,\dots , X_n$ be independent random variables, such that $X_i$ takes his values in $[a_i,b_i]$ almost surely for all $i \leq n$. Then for any positive x,
  \begin{eqnarray*}
  \prob{}{\sum_{i=1}^n \Big(X_i - \esp{X_i}\Big) \geq x} & \leq & \exp\left( - \frac{2
x^2}{\sum_{i=1}^n(b_i -  a_i)^2}  \right).
  \end{eqnarray*}
\end{lmm}
We may also apply this inequality to $-X_i$ so as to get a very probable lower bound on the sum of $X_i$.

This is enough to prove:
\begin{thrm}
\label{deter}
  Let $\rho$ be a density operator. Assume that each $f_i$ is
bounded,  where $(f_i)_{i \in \mathcal{I}}$ is the dual basis of $(e_i)_{i \in \mathcal{I}}$,
as defined at the beginning of this section. Let $M_i = \sup_{(x,\phi) \in
\mathbb{R}\times[0,\pi]} f_i(x,\phi) - \inf_{(x,\phi) \in \mathbb{R}\times[0,\pi]}
f_i(x,\phi)$. Let $(x_i)_{i\in\mathcal{I}}$ be a family of positive real numbers such that
$\sum_{i\in\mathcal{I}} \exp(-x_i) = \sigma < \infty$. Let
\begin{eqnarray}
\label{penalty}
\pen_n(m) & = & \sum_{i \in \mathcal{I}_m} (1+\epsilon) \left( \ln(M_i) + \frac{x_i}{2}\right) \frac{M_i^2}{n}.
\end{eqnarray} 
  Then the penalized projection estimator satisfies:
\begin{eqnarray}
\label{eqdeter}
\esp{\frac{\epsilon}{2+\epsilon}\norm{\hat{\rho}^{(n)} - \rho}^2} & \leq & \inf_{m\in \mathcal{M}}\left(1 +
\frac{2}{\epsilon}\right) d^2(\rho,m) + 2 \pen_n(m) + \frac{(1+\epsilon)\sigma}{n}.
\end{eqnarray}
  
\end{thrm}

{\bf Remark:} Here the penalty depends only on the subspace spanned by the
model $m$. So it is the same whether $\mathcal{M}$ is small or large. The best we
can do is then to take $\mathcal{M}= \mathcal{P}(\mathcal{I})$, that is to choose for every vector
$e_i$ whether to keep the estimated coordinate $\hat{\rho}_i$ or to put it to
zero. In other words we get a hard-thresholding estimator:
\begin{eqnarray*}
\hat{\rho}^{(n)}& = &\sum_{i\in\mathcal{I}} \hat{\rho}_i \bf{1}_{|\hat{\rho}_i|>\alpha_i} e_i
\end{eqnarray*} 
with
\begin{eqnarray}
\label{seuil_hoeffding}
\alpha_i& =& \sqrt{(1+\epsilon)\left(\ln(M_i) + \frac{x_i}{2}\right)}\frac{M_i}{\sqrt{n}}
\end{eqnarray}

\begin{proof}

  Considering (\ref{esperancereelle2}), we have only to bound appropriately
$\esp{\sup_m \left((1+\epsilon)\chi_n^2(m) - \pen(m)\right)\vee 0} $. 

Now, by (\ref{chi2}) and (\ref{penalty}), both $\chi_n^2(m)$ and $\pen_m$ are a
sum of terms over $m$. As the positive part of a sum is smaller than
the sum of the positive parts, we obtain:
\begin{multline*} 
 \esp{\sup_m \left\{(1+\epsilon)\chi_n^2(m) - \pen(m)\right\}\vee 0} \\  
 \begin{aligned}
&   \leq   
\esp{\sup_m \left\{\sum_{i\in m}\left((1+
\epsilon)\left(\hat{\rho}_i -
\rho_i \right)^2 - \alpha_i^2 \right\}\vee 0\right)} \\
&= 
\sum_{i\in\mathcal{I}}\esp{\left\{(1+\epsilon)\left(\frac1{n}\sum_{k=1}^{n}f_i(x_k,\phi_k)  -
\rho_i \right)^2 - (1+\epsilon) \left( \ln(M_i) +
\frac{x_i}{2}\right) \frac{M_i^2}{n}
 \right\}\vee 0}. 
\end{aligned}
\end{multline*}
Each of the expectations is evaluated using the following formula, valid for any positive function~$f$:
\begin{eqnarray}
\label{form_esp}
\esp{f} = \int_0^{\infty} \prob{}{f(x) \geq y} \dd y.
\end{eqnarray}
Remembering \eqref{seuil_hoeffding} we notice that the inequality
\begin{eqnarray*}
\left\{(1+\epsilon)\left(\frac1{n}\sum_{k=1}^{n}f_i(x_k,\phi_k)  -
\rho_i \right)^2 - (1+\epsilon) \left( \ln(M_i) +
\frac{x_i}{2}\right) \frac{M_i^2}{n}
 \right\}\vee 0 & \geq & y
\end{eqnarray*}
is equivalent to 
\begin{eqnarray*}
\left|\frac1{n}\sum_{k=1}^{n}f_i(x_k,\phi_k)  -
\rho_i \right| & \geq & \sqrt{\frac{\alpha_i^2 + y}{1+ \epsilon}}.
\end{eqnarray*}
We may then conclude, using Hoeffding's inequality on the second line and the value \eqref{seuil_hoeffding} of $\alpha_i$ on the fourth line:
\begin{eqnarray*}
 \esp{\sup_m \left\{(1+\epsilon)\chi_n^2(m) - \pen(m)\right\}\vee 0} 
& \leq &   \sum_{i\in\mathcal{I}} \int_{0}^{\infty} \prob{}{\left|\frac1{n}\sum_{k=1}^{n}f_i(x_k,\phi_k)  - \rho_i \right| \geq \sqrt{\frac{\alpha_i^2 + y}{1+ \epsilon}}}  dy \\
 & =   &  \sum_{i\in\mathcal{I}} \int_0^{\infty} 2 \exp\left(- \frac{2 n (\alpha_i^2 + y)}{(1+\epsilon) M_i^2}\right) \dd y \\
& = & \sum_{i\in \mathcal{I}} 2 \exp\left(- \frac{2 n \alpha_i^2}{(1+\epsilon) M_i^2}\right) \frac{(1+\epsilon) M_i^2}{2 n} \\
& = & \frac{1 + \epsilon}{n} \sum_{i\in\mathcal{I}} \exp(- x_i) \\
& = & \frac{(1 + \epsilon) \sigma}{n}. 
\end{eqnarray*}

\end{proof}

\subsection{Random penalty}
\label{ranpen}

The most obvious way to improve on Theorem \ref{deter} is to use sharper
inequalities than Hoeffding's. Indeed the range of $f_i$ might be much larger
than its standard deviation, so that we gain much by using Bernstein's inequality:
\begin{lmm}{\bf: Bernstein's inequality}
Let $X_1,\dots,X_n$ be independent, bounded, random variables. Then with
\begin{align*}
M & =\sup_i\norm{X_i}_{\infty}, & v&=\sum_{i=1}^n\esp{X_i^2},
\end{align*}
for any positive x
\begin{eqnarray*}
\prob{}{\sum_{i=1}^n (X_i -  \esp{X_i}) \geq \sqrt{2vx}+\frac{M}{3}x} & \leq & \exp(-x). 
\end{eqnarray*}  
\end{lmm}
With this tool, we may devise a hard-thresholding estimator where the thresholds are data-dependent:
\begin{thrm}
\label{alea}
  Let $(y_i)_{i\in \mathcal{I}}$ be positive numbers such that $\sum_{i\in\mathcal{I}} e^{-y_i} = \sigma < \infty$. Let then 
\begin{eqnarray*}
  x_i & = & 2\ln(\norm{f_i}_{\infty}) + y_i.
\end{eqnarray*}
Let the penalty be a sum of penalties over the vectors we admit in the model.
That is, for any $\delta \in (0,1)$, for any $i\in \mathcal{I}$, define
\begin{eqnarray}
\label{pen_i}
\pen_n^i & = & \frac{1+\epsilon}{n} \left( \sqrt{\frac{2}{1-\delta}x_i \left(\prob{n}{f_i^2}+\frac1{n} \norm{f_i}^2_{\infty}(\frac1{3}+\frac1{\delta})x_i\right)} + \frac{\norm{f_i}_{\infty}}{3\sqrt{n}}x_i \right)^2
\end{eqnarray}
and the penalty of the model $m$:
\begin{eqnarray*}
\pen_n(m) & = &\sum_{i\in m} \pen^i_n.
\end{eqnarray*}

Then there is a constant C such that:
\begin{eqnarray*}
\esp{\left(\frac{\epsilon}{2+\epsilon}\norm{\hat{\rho}^{(n)} - \rho}^2 -  \left( \inf_{m\in \mathcal{M}_n}\left(1 + \frac{2}{\epsilon}\right) d^2(\rho,m) + 2 \pen_n(m) \right)\right)\vee 0} & \leq & \frac{C\sigma}{n}
\end{eqnarray*}
where $\mathcal{M}_n$ is the set of models $m$ for which $i\in m \rightarrow x_i\leq n$.   
  
\end{thrm}

{\bf Remark:} As with the deterministic penalty, we end up with a
hard-thresholding estimator. Morally, that is, forgetting all the small $\delta$ whose origin is technical, the threshold is
\begin{eqnarray*}
 \sqrt{\frac{2\prob{n}{f_i^2}\ln\norm{f_i}_{\infty}^2}{n}}.
\end{eqnarray*}

\begin{proof}
Once again we have to dominate the right-hand side of (\ref{esperancereelle2}).
We first subtract and add, inside that expression, what could be seen as a target for the penalty. Writing 
\begin{align}
\label{not1}
M_i & = \norm{f_i}_{\infty}, & v_i & = n \esp{f_i^2}, & \alpha_i & =  \sqrt{2v_i x_i}  + \frac{M_i}{3}x_i
\end{align}
we have 
\begin{multline}
\label{step1}
 \esp{\sup_m \left((1+\epsilon)\chi_n^2(m) - \pen(m)\right)\vee 0}   \\  
  \leq  \esp{\sup_m (1+\epsilon)\left(\chi_n^2(m) -
\sum_{i\in m}\frac{1}{n^2}
\alpha_i^2\right)\vee 0}  
 +\, \esp{\left(\sum_{i\in m}\frac{1+\epsilon}{n^2}
\alpha_i^2- \pen(m)\right)\vee 0}.
\end{multline}

Using \eqref{chi2}, we bound the first term as a sum of expectations.
\begin{eqnarray*}
 \esp{\sup_m (1+\epsilon)\left(\chi_n^2(m) -
\sum_{i\in m}\frac{1}{n^2}
\alpha_i^2\right)\vee 0} &
 \leq & (1 + \epsilon) \sum_{i\in m} \esp{ \left(\left|\rho_i-\frac1n\sum_{k=1}^n
f_i(x_k,\phi_k)       \right|^2 -  \frac{1}{n^2}
\alpha_i^2\right)\vee 0}. 
\end{eqnarray*}

We now bound each of these expectations using \eqref{form_esp}.
\begin{eqnarray} 
\esp{ \left(\left|\rho_i-\frac1n\sum_{k=1}^n
f_i(x_k,\phi_k)       \right|^2 -  \frac{1}{n^2}
\alpha_i^2\right)\vee 0}  
\label{avec_prob}
& = & \int_{0}^{\infty}  \prob{}{\left|\rho_i-\frac1n\sum_{k=1}^n
f_i(x_k,\phi_k)       \right| \geq \sqrt{y + \frac{\alpha_i^2}{n^2}}} \dd y.
\end{eqnarray}
We change variables in the integral, choosing $\xi$ defined by:
\begin{eqnarray}
\label{change_var}
\sqrt{y + \frac{\alpha_i^2}{n^2}} & = & \frac{\sqrt{2v_i \xi} + \frac{M_i}3 \xi}{n^2}.
\end{eqnarray}
Using Bernstein's inequality, the integrand in \eqref{avec_prob} is upper bounded by $2\exp(-\xi)$. Given the value of $\alpha_i$ \eqref{not1}, the range of the integral is now from $x_i$ to $\infty$. Finally, taking the square on both sides of \eqref{change_var}, then using \eqref{abaabb}, we get:
\begin{eqnarray*}
\dd y & = & 2 \frac{\sqrt{2v_i \xi} + \frac{M_i}3 \xi}{n^2} \left(\frac{M_i}{3} + \frac{\sqrt{2v_i}}{2\sqrt{\xi}}\right)\dd \xi \\
 & = & \frac{2}{n^2}\left(v_i + \frac{M_i^2}{9}\xi + \frac{M_i}{2} \sqrt{2v_i} \sqrt{x}\right) \dd \xi \\
& \leq & \frac{2}{n^2}\left(2 v_i + \frac{11 M_i^2}{18} \xi \right) \dd \xi.
\end{eqnarray*}

Hence
\begin{eqnarray}
\esp{ \left(\left|\rho_i-\frac1n\sum_{k=1}^n
f_i(x_k,\phi_k)       \right|^2 -  \frac{1}{n^2}
\alpha_i^2\right)\vee 0} 
& \leq & \frac{4}{n^2} \int_{x_i}^{\infty} \exp(-\xi)  \left(2 v_i + \frac{11 M_i^2}{18} \xi \right)  \dd \xi \notag \\
\label{terme_1}
& = & \frac{4}{n^2}\left(2 v_i + \frac{11 M_i^2}{18}(1 + x_i) \right) \exp(-x_i).
\end{eqnarray}

Let us now look over the second term of \eqref{step1}. We notice, through \eqref{pen_i} and \eqref{not1}, that this term is of the form:
\begin{eqnarray*}
 \frac{1 + \epsilon}{n^2}\sum_{i\in m} \esp{\left(\left(a_i + \frac{M_i x_i}{3}\right)^2 -\left(b_i + \frac{M_i x_i}{3}\right)^2 \right) \vee 0} 
& \leq & \frac{1 + \epsilon}{n^2}\sum_{i\in m} \esp{2 \left( a_i^2 - b_i^2 \right) \vee 0},
\end{eqnarray*} 
with 
\begin{eqnarray*}
a_i^2 - b_i^2 & = & 2 v_i x_i - \frac{2}{1 - \delta}\left(n \prob{n}{f_i^2} x_i +M_i^2\left(\frac1{3}+\frac1{\delta}\right) x_i^2 \right).
\end{eqnarray*}
Using again \eqref{form_esp}, we end up with:
\begin{multline}
\esp{\left(\sum_{i\in m}\frac{1+\epsilon}{n^2} \alpha_i^2- \pen(m)\right)\vee 0} \\
\label{avant_integ}
 \leq  \frac{1 + \epsilon}{n^2} \sum_{i\in m} \frac{2}{1 - \delta}x_i \int_{0}^\infty \prob{}{(1 - \delta) v_i -\left(n \prob{n}{f_i^2}  +M_i^2\left(\frac1{3}+\frac1{\delta}\right) x_i \right) \geq y} \dd y.
\end{multline}

We can again make use of Bernstein's inequality:
\begin{eqnarray*}
\prob{}{v_i - \sum_{k=1}^n f_i^2(x_k,\phi_k)\geq \sqrt{2 n
\esp{f_i^4}\xi}+\frac{\norm{f_i^2}_{\infty}\xi}{3}} & \leq & \exp(-\xi).
\end{eqnarray*}
Noticing that $f_i^2$ is non-negative everywhere, so that $\esp{f_i^4}\leq\esp{f_i^2}\norm{f_i^2}_{\infty}$, and using (\ref{abaabb}), we get:
\begin{eqnarray*}
  \prob{}{(1-\delta) v_i \geq 
n\prob{n}{f_i^2}+M_i^2\left(\frac1{3}+\frac1{\delta}\right)\xi} & \leq & \exp(-\xi).
\end{eqnarray*}

Recalling \eqref{avant_integ}, we get
\begin{eqnarray*}
 \int_{0}^\infty \prob{}{(1 - \delta) v_i -\left(n \prob{n}{f_i^2}  +M_i^2\left(\frac1{3}+\frac1{\delta}\right) x_i \right) \geq y} \dd y 
& = &
\int_{0}^\infty  \exp\left(- x_i - \frac{y}{M_i^2\left(\frac1{3}+\frac1{\delta}\right)}\right) \dd y \\
& = & \exp(- x_i) M_i^2\left(\frac1{3}+\frac1{\delta}\right)  \exp\left(-\frac{x_i}{M_i^2\left(\frac1{3}+\frac1{\delta}\right)}\right) \\
& \leq & \exp(-y_i) \left(\frac1{3}+\frac1{\delta}\right). 
\end{eqnarray*}

With that and (\ref{terme_1}), we are left with:
\begin{eqnarray*}
  \esp{\sup_m \left\{(1+\epsilon)\chi_n^2(m) - \pen(m)\right\}\vee 0} 
 & \leq & \frac{C}{n^2}\sum_{i\in\mathcal{I}} e^{-x_i}(v_i+ M_i^2(1+x_i)) + x_i e^{-y_i}.
\end{eqnarray*}

Replacing $x_i$ with its value, and overestimating $v_i$ by $n M_i^2$ we obtain (under the condition that $2\ln M_i + y_i \leq n$):
\begin{eqnarray*}
 \esp{\sup_m \left\{(1+\epsilon)\chi_n^2(m) - \pen(m)\right\}\vee 0}   
 & \leq & C \left(\frac{\sigma}{n}+ \frac{\sigma}{n^2}\right).
\end{eqnarray*}

\end{proof}

{\bf Remark:} The logarithmic factor in the penalty (that would not be here if we took only the variance) comes from the multitude of models allowed by a hard-thresholding estimator. By selecting fewer models (for example the square matrices obtained by truncating the density operator) and using a random penalty, we can get rid of this term. However, crafting the penalty requires much more work and more powerful inequalities (Talagrand's). An interested reader may have a look at the section 3.4 of \cite{DEA}.

\subsection{Applications with two bases}

  We now give two bases that are reasonable when applying these theorems. As can be
seen from (\ref{oracle}), a good basis should approximate well any density operator (so that the bias term gets low fast when $m$ is big), with dual vectors having a low variance. With the first of the two bases, we have this interesting phenomenon that we obtain a polynomial convergence rate under the mere physical hypothesis that the state has finite energy.

\subsubsection{Photon basis}
\label{photon}
  
  Here we shall take as our $(e_i)_{i\in \mathcal{I}}$ a slight variation of the matrix entries of our density operator with respect to the Fock basis \eqref{hermite}.

  More precisely, we worked in the previous subsections with real  coefficients. To apply Theorems \ref{alea} and \ref{deter}, we then need to parametrize $\rho$ with real coefficients. The matrix entries are \emph{a priori} complex. However, using the fact that $\rho$ is self-adjoint, we may separate the real and imaginary parts.  

We use a double index for $i$ and define the orthonormal basis, denoting by $E_{j,k}$ the null matrix except for a $1$ in case $(j,k)$:
\begin{eqnarray*}
e_{j,k} & = & \left\{ \begin{array}{ll}
              \frac1{\sqrt 2} (E_{j,k} + E_{k,j}) & \quad \textrm{if } j<k \\
              \frac{i}{\sqrt 2} (E_{k,j} - E_{j,k}) & \quad \textrm{if } k<j \\
              E_{j,j} & \quad \textrm{if} j=k
              \end{array} \right. .
\end{eqnarray*}

Then, using a tilde to distinguish it from the matrix entries, with $\tilde{\rho}_{j,k}  = \langle \rho, e_{j,k} \rangle$,we have 
\begin{eqnarray*}
 \langle \psi_j , \rho \psi_k \rangle  = \left\{ \begin{array}{ll} 
                                                    \frac{1}{\sqrt{2}}(\tilde{\rho}_{j,k} + i\tilde{ \rho}_{k,j}) & \mathrm{if} \,\, j<k \\
                                                    \frac1{\sqrt{2}}(\tilde{\rho}_{k,j} - i \tilde{\rho}_{j,k}) & \mathrm{if}  \,\, j>k \\
						    \tilde{\rho}_{j,j}
& \mathrm{if}  \,\, j=k. \\
                                                   \end{array}
						\right.                                  
\end{eqnarray*}

The associated $\tilde{f}_{j,k}$ are well-known. They are a slight variation of the usual ``pattern functions'' (see Appendix \ref{QHD}, and \eqref{pattern} therein), the behaviour of which may be found in \cite{Gu.Gill.}. Notably, we know that:
\begin{equation}
\label{norm}
\sum_{j,k=0}^{N} \norm{f_{j,k}}^2_{\infty} \leq CN^{7/\!3}.
\end{equation}
As the upper bounds on the supremum of $\tilde{f}_{j,k}$ may not be sharp,
the best way to apply the above theorems (especially Theorem (\ref{deter})) would probably be to tabulate these maxima for the $(j,k)$ we plan to use.

 The interest of this basis is that it is a priori adapted to the structure of our problem: if we have a bound on the energy (let's say it is lower than $H+\frac1{2}$), we get worst-case estimates on the convergence speed with the deterministic penalty: indeed, the energy of a state $\rho$ may be written $\frac1{2} + \sum_j j\rho_{j,j}$, so that
\begin{eqnarray*}
\sum_{j\geq N} \tilde{\rho}_{j,j} & \leq & \frac{H}{N}.
\end{eqnarray*}
  Moreover, by positivity of the operator, 
\begin{eqnarray*}
\tilde{\rho}_{j,k}^2+\tilde{\rho}_{k,j}^2 & \leq &
\tilde{\rho}_{j,j}\tilde{\rho}_{k,k}.
\end{eqnarray*}
If we look at the models $N$ such that $\mathcal{I}_N = \{(j,k): j<N, k<N \}$, we can get:
\begin{eqnarray*}
d^2(\rho,N) & \leq &  \sum_{j,k=0}^{\infty} \tilde{\rho}_{j,k}^2 - \sum_{j,k=0}^{N} \tilde{\rho}_{j,k}^2  \\
                 & \leq &  (\sum_{j\geq 0} \tilde{\rho}_{j,j})^2 - (\sum_{j=0}^{N} \tilde{\rho}_{j,j})^2  \\
		 & \leq &  1 - (1-\frac{H}{N})^2 \\
		 & \leq &  \frac{2H}{N}
\end{eqnarray*}
where we have used that the density operator has trace one.

We substitute in (\ref{eqdeter}) and get:
\begin{eqnarray*}
  \esp{\norm{\hat{\rho}^{(n)} - \rho}^2} & \leq & C\left(\frac{H}{N} +  \pen_n(N) +
\frac1{n}\right). 
\end{eqnarray*}

Now, using the bounds on infinite norms (\ref{norm}), the penalty is less than:
\begin{eqnarray*}
\pen_n(N) & = & C \frac{N^{7\!/3}\ln(N)}{n}.
\end{eqnarray*}
Optimizing in $N$ ($N=C(Hn)^{3\!/10}$), we get
\begin{eqnarray}
\label{vit_app}
  \esp{\norm{\hat{\rho}^{(n)} - \rho}^2}  & \leq & C H^{7\!/10}\ln(H)n^{-3\!/10}\ln(n). 
\end{eqnarray}

This estimate  holds true for any state and is non-asymptotic. It is
generally rather pessimistic, though. For many classical states, such as
squeezed states or thermal states, $\rho_{j,j}\equiv A \exp(-B/\!n)$, the same
calculation  yields a rate for the square of the $L^2$-distance as
$n^{-1}\ln(n)^\beta$ for some $\beta$. In such a case, the penalized estimator 
automatically converges at this latter rate.

\subsubsection{Wavelets}

  Another try could be to use functions known for their good approximations properties. To this end we look at the Wigner function and write it in a wavelet basis.
  
  Recall that wavelets on $\mathbb{R}$ are an orthonormal basis such that all functions
are  scaled translations of a same function, the mother wavelet. In multiscale
analysis, we use a countable basis  $\psi_{j,k}: x \mapsto 2^{j\!/2}\psi_{0,0}(2^j x +k)$, for j and k integers. Let $\mathcal{V}_i=\{\psi_{j,k} : j\leq i\}$. There is a $\phi$, called father wavelet, such that the $\phi_{k}(x)=\phi(x+k)$ for $k\in  \mathbb{Z}$ are a basis of the vector space generated by all the wavelets of  larger  or equal scale, that is $\mathcal{V}_0$. We may choose them with compact support, or localized both in frequency and position. So they harvest local information and can fetch this whatever the regularity of the function to be approximated, as they exist at several scales.

  From a one-dimensional wavelet basis $\psi_{j,k}: x \mapsto 2^{j\!/2}\psi_{0,0}(2^j x +k)$, $C^3$ and zero mean, with a father wavelet $\phi_{j,k}$, also $C^3$, we shall make a tensor product basis on $L^2(\mathbb{R}^2)$: let  $I=(j,k,\epsilon)$ be indices, with j integer (scale), $k = (k_x, k_y) \in  \mathbb{Z}^2$ (position), and $\epsilon\in{0,1,2,3}$. Let
\begin{eqnarray*}
  \Psi_{I}(x,y) = \left\{
                        \begin{array}{ll}
			\phi_{j,k}(x)\phi_{j,k}(y) & \mathrm{if } \ \ \epsilon=0 \\
			\phi_{j,k}(x)\psi_{j,k}(y) & \mathrm{if } \  \ \epsilon=1 \\
			\psi_{j,k}(x)\phi_{j,k}(y) & \mathrm{if  } \ \  \epsilon=2 \\
			\psi_{j,k}(x)\psi_{j,k}(y) & \mathrm{if  } \ \  \epsilon=3 \\
			\end{array}
                   \right.
\end{eqnarray*}
  
We may then define a multiscale analysis from the one-dimensional one (written  $\mathcal{V},\mathcal{W}$): $V_0 =\overline{\mathcal{V}_0\otimes\mathcal{V}_0}$ and for all $j \in  \mathbb{Z}$, $V_{j+1}=V_j \oplus W_j$, so that $W_{j+1} = \overline{\mathcal{V}_j\otimes\mathcal{W}_j}\oplus \overline{\mathcal{W}_j\otimes\mathcal{V}_j} \oplus \overline{\mathcal{V}_j\otimes\mathcal{W}_j}$.

For any j, $V_j \cup \bigcup_{k\geq j}W_k$ is then an orthonormal basis of $L^2(\mathbb{R}^2)$. We  hereafter choose our models as subspaces spanned by finite subsets of the basis vectors for well-chosen j.

It can be shown that: 
\begin{eqnarray*} 
\gamma_I(x,\phi) &  = & \frac1{4\pi} \int_{-\infty}^\infty \left|u\right| \hat{\Psi}_I(u \cos \phi, u \sin \phi) e^{ixu} du
\end{eqnarray*}
fulfills this property: 
\begin{eqnarray*} 
[\gamma_I,Kf] & = & \ps{\Psi_I}{f}. 
\end{eqnarray*}
Noticing that
\begin{eqnarray*} 
\gamma_I(x,\phi) & = & 2^j\gamma_{0,0,\epsilon}(2^j x - k_x \cos \phi - k_y \sin \phi,
\phi),
\end{eqnarray*}
we see that these functions have the same dilation properties as wavelets, and they are ``translated'' in a way that depends on $\phi$, through sinusoids. Their normalizations, though, explode with $j$; this derives from inverting the Radon transform being an ill-posed problem.

We can now apply Theorem \ref{alea}. Before doing so, though, we  restrict ourselves to a finite subdomain of $\mathbb{R}^2$, which we  denote $\mathcal{D}$, and put the Wigner function to zero outside this domain, that we should choose big enough to ensure this does not cost too much.
 
 Then, $\mathcal{M}$ is the set of all models characterized by 
\begin{eqnarray*} 
m & = & \left\{ (j_1,k,0): 2^{j_1}k\in \mathcal{D} \right\} \cup \left\{(j,k,\epsilon): (j,k,\epsilon)\in \mathcal{I}_{m}' \subset \{ (j,k,\epsilon): \epsilon=1;2;3 , j_1<j<j_0 , 2^j k\in\mathcal{D} \}\right\}.
\end{eqnarray*}
To have good approximating properties, we choose $2^{j_1}\equiv n^{1\!/7}$ and $2^{j_0}\equiv \frac{n}{(\ln n)^2}$. The projection estimator within a model is then:
\begin{eqnarray*}
\hat{f} & = & \sum_{I\in m} \alpha_I \Psi_I 
\end{eqnarray*}
with 
\begin{eqnarray*}
\alpha_I & = & \frac1{n}\sum_{i=1}^n \gamma_I(x_i,\phi_i).
\end{eqnarray*}

Denoting $B_{\epsilon}= \norm{\gamma_{0,0,\epsilon}}_\infty$, the translation of Theorem $\ref{alea}$ gives (notice that applying (\ref{deter}) would be awkward, as the variance of $\gamma_I$ is like $2^j$ whereas its maximum is like $2^{2j}$):

\begin{thrm}

  Let $y_I$ be such that $\sum_{I}\exp(-y_I) = \sigma \leq \infty$. For example $y_I=j$. Let then:
\begin{eqnarray*}
  x_I & = &  2(j + \ln(B_{\epsilon})) + y_{I}. 
\end{eqnarray*}
  We choose an $\alpha \in (0,1)$ and the penalty (and restraining ourselves to the $m$ such that $I\in m \rightarrow x_I\leq n$):
\begin{eqnarray*}
  \pen(m) & = & \frac{1+\epsilon'}{n}\sum_{I\in \mathcal{M}} 2\left(
\sqrt{\frac{2}{1-\alpha}x_{I}\left(\prob{n}{\gamma_{I}^2}+\frac1{n}
2^{2j}B_{\epsilon}^2\Big(\frac1{3}+\frac1{\alpha}\Big)x_{I}\right)} +
\frac{2^{j}B_{\epsilon}}{3\sqrt{n}}x_{I} \right)^2.
\end{eqnarray*}
  
 Then there is a constant C such that:
\begin{eqnarray}
\label{ecafi}
\esp{\left\{\frac{\epsilon}{2+\epsilon}\norm{\rho-\hat{\rho}^{(n)}}^2 - \left (\inf_{m \in \mathcal{M}}
\left( 1 + \frac{2}{\epsilon}\right) d^2(\rho,m) + 2 \pen_n(m)\right)\right\}\vee 0}
& \leq & \frac{C \sigma}{n} + C\frac1{n}2^{2j_1}. 
\end{eqnarray}

\end{thrm}

\begin{proof}

First it's easily checked that $x_I = 2 \ln(\norm{\gamma_I}_{\infty}) + y_I$. Second $\sum_{I}\exp(-j) = C \sum_{j}2^j \exp(-j) <\infty$ implies that $y_I=j$ does indeed the work, as there are at most $C 2^j$ wavelets at scale j whose support meet $\mathcal{D}$.

The last term is the variance of $ \hat{a}_{j_1,k,0} $, corresponding to the vectors that are in every model.:
\begin{eqnarray*}
\frac1{n}\mathbb{V}\left[\sum_{2^{j_1}k \in \mathcal{D} }\gamma_{j_1,k,0}\right] & \leq & \frac1{n}\esp{\sum_{2^{j_1}k \in \mathcal{D} }\gamma_{j_1,k,0}^2} \\
                                                             & \leq & \frac1{n}\sum_{2^{j_1}k\in \mathcal{D}}\int_{\mathbb{R}\times [0,\pi]}\gamma_{j_1,k,0}^2(x,\phi) dx \frac{d\phi}{\pi} p_{\rho}(x,\phi)\\
							     & =    & \frac1{n}\sum_{2^{j_1}k\in \mathcal{D}}\int_{\mathbb{R}}\gamma_{j_1,k,0}^2(x,0)\int_0^{\pi} p_{\rho}(x-k_x \cos \phi -k_y \sin \phi,\phi) dx \frac{d\phi}{\pi} \\
							     & =    & C\frac1{n}2^{2j_1} 
\end{eqnarray*}
where we have used that for all x and k, 
$
\int_0^{\pi} p_{\rho}(x-k_x \cos \phi -k_y \sin \phi,\phi) \frac{d\phi}{\pi}
$ 
is less than a constant about $1.086$. Indeed, the translation of a Wigner function is still the Wigner function of a state, so that we may take $k=0$. Then 
\begin{eqnarray*}
\int_0^{\pi} p_{\rho}(x-k_x \cos \phi -k_y \sin \phi,\phi) \frac{d\phi}{\pi} & \leq & \sup_{i,x} |\psi_i(x)|^2
\end{eqnarray*}
and the upper bound on this supremum is due to Cramér (10.18.19 in \cite{Erdelyi}).
\end{proof}

\emph{Remarks:} As the variance of $\gamma_I$ goes like $2^{j}$ the threshold might be seen as $C 2^{j\!/2}\sqrt{\frac{j}{n}}$. This is the estimator studied in \cite{Cavalier2}, for a general Radon transform (i.e. not a Wigner function). 



 The role of the approximation speed is apparent in (\ref{ecafi}). Articles
like \cite{Cavalier2} show that this strategy is asymptotically (quasi)-optimal
for approximating a function in a Besov ball. However, this is no proof of the
efficiency in our case, as the set of Wigner functions is not a Besov ball.
There is still some work in approximation theory needed there. In particular, we do not know if a statement similar to \eqref{vit_app} can be proven.
 
\bigskip

Finally, notice that we may combine projection estimators: as the contrast function
is the same for any basis we are working with, keeping the same penalizations,
we could find an estimator that is almost the best among those built with the photon basis and those with the wavelet basis simultaneously (just add a $\ln(2)$ to $\sigma$). In other words, we do not have to choose beforehand which basis we use. Moreover an estimator allowing for the two bases would satisfy \eqref{vit_app}

\subsection{Noisy observations}
\label{noise}

  The situation we have studied was very idealized: we did not take any noise into account. In practice, a number of photons fail to be detected. These losses may be quantified by one single coefficient $\eta$ between $0$ (no detection) and $1$ (ideal case). We suppose it to be known. 

There are several methods to recover the state from noisy observations. One
consists in calculating the density matrix as if there was no noise, and then
apply the Bernoulli transformation with factor $\eta^{-1}$. We can also use
modified pattern functions \cite{D'Ariano.2}. Or we can approximate the Wigner function with a kernel estimator that performs both the inverse Radon transform and the deconvolution \cite{But.Gu.Art.}. The former two methods fail if the observations are too noisy ($\eta \leq \frac 12$), but the latter is asymptotically optimal for all $\eta$ over wide classes of Wigner functions.

   This noise can be seen as a convolution of the result $(X,\Phi)$ with a gaussian of variance depending on $\eta$:
\begin{eqnarray*}
  p_{\rho}^{\eta}(y,\phi) & = & \frac{1}{\sqrt{\pi(1-\eta)}}\int_{-\infty}^{\infty} p_{\rho}(x,\phi)\exp\left(-\frac{\eta}{1-\eta}\left(x-\eta^{-1/\!2}y\right)^2\right)dx
\end{eqnarray*}
or equivalently in terms of generating functions
\begin{eqnarray*}
\int p_\rho^\eta (x,\phi)e^{irx}dx & = & e^{-\frac{1-\eta}{8\eta}r^2}\int
p_\rho(x,\phi) e^{irx}dx. 
\end{eqnarray*}

We can use the methods described above and then use the Bernoulli transform. For free, we may also use the modified pattern functions $f_{j,k}^{\eta}$ knowing $f_{j,k}$. Explicitly we see that the dual basis of the matrix entry $\rho_{j,k}$  becomes:
\begin{eqnarray*}
f_{j,k}^{\eta}(x,\phi) & = & \frac 1 {2\pi} \int dr e^{\frac{1-\eta}{8\eta}r^2}
\int dy f_{j,k}(y,\phi) e^{iry}. 
\end{eqnarray*}
The reason why one needs $\eta>\frac 12$ is for this Fourier transform to be well defined. 

And we can again apply Theorems \ref{deter} and \ref{alea} with the dual basis $ \tilde{f}_{j,k}^{\eta}$.  

  Obtaining results with high noise $\eta\leq \frac 12$ is harder. We would
need to introduce a cut-off $h$ within the inverse Fourier transform (and
therefore a bias). Using the same $h$ as in \cite{But.Gu.Art.} would ensure this bias
$b(\rho,h)$ is asymptotically reasonable. We could then reuse Theorems
\ref{deter} and \ref{alea} to have an ``almost best'' approximation of $\rho +
b(\rho,h)$ within a set of models, for finite samples. Careful examination would
then be required to check the variance (or the penalties) go to $0$ as $n$ and
$h(n)$ go to infinity. Moreover, we would need to translate conditions on the
Wigner function into conditions on the density operator to see whether we can reproduce the asymptotic optimality results of Butucea {\em et al.} with model selection in the Fock basis (or any other basis chosen and studied {\em a priori}).

\section{Maximum likelihood estimator}
\label{EMV}

  Projection estimators are not devoid of defects: the variance of empirical coefficients might be high, and the convergence therefore rather slow, the estimator is not a true density matrix... Especially, the trace is probably not one, though this could be fixed easily enough. We can diagonalize the estimated density matrix, replace  the negative eigenvalues with $0$,  and divide by the trace.
  
  Anyhow, there are other types of estimator that automatically yield density matrices. One such estimator is the maximum likelihood estimator, which selects the nearest point of the empirical probability measure in a given model for the Kullback-Leibler distance (which is not a true distance as it is not symmetric). Recall that the Kullback-Leibler distance between two probability measures is:
\begin{eqnarray*}
 K(p,q) & = & \int \ln\left(\frac{p(x)}{q(x)}\right)p(x) dx.
\end{eqnarray*}
In other words, the maximum likelihood estimator is
\begin{eqnarray*}
 \underset{\tau \in \mathcal{Q}}{\arg\min} \sum_{l=1}^{n} - \ln
p_{\tau}(X_l,\Phi_l)
\end{eqnarray*}
where $\mathcal{Q}$ is any set of density operators (such that the minimum exists). This way, it is automatically a true density operator. A practical drawback is that calculating it is very power-consuming.
 
As $\gamma_n(\cdot)\rightarrow - \int \ln(p_{\cdot}) d_{p_\rho}$, we have defined a
minimum contrast estimator in the sense of section \ref{31}. Much like for
projection estimators, the Kullback distance thus estimated might be overly
optimistic, and all the more when $\mathcal{Q}$ is big. Indeed, if $ \mathcal{Q}$ is the set of all density operators, then there is no minimizer of the distance with the empirical
distribution; however when we take only finite-dimensional models, such as 
\begin{eqnarray}
\label{tronqueN}
 \mathcal{Q}(N) & = & \left\{ \tau \in \mathcal{S}(L^2(\mathbb{R})) : \tau_{j,k} = 0 \,\,\,\mathrm{for ~ all }\,
j > N\,\,\mathrm{or} \,k > N \right\}, 
 \end{eqnarray}
then the minimum is attained by compactness. Here the matrix entries $\tau_{j,k}$ are taken in the Fock basis (\ref{hermite}).

We then have to define a penalty for choosing (almost) the best model. To do so, we  make use of a (slightly simplified but sufficient for our needs) version of a theorem by Massart \cite{Mas}, but we need a few definitions before stating it.

First we  need a distance with which to express our results, and it is not the Kullback-Leibler, but the Hellinger distance.
The Hellinger distance between two probability densities  is defined in relation with the $L^2$-distance of the square roots of these densities:
\begin{eqnarray}
\label{defhell}
  h^2(p,q) & = & \frac1{2} \int \left(\sqrt{p}-\sqrt{q}\right)^2.
\end{eqnarray}
This distance does not depend on the underlying measure. The following relations are well-known:
\begin{eqnarray}
\label{comparehellinger}
  \frac1{8}\|p-q\|_1^2 & \leq &  h^2(p,q)  \leq  \frac1{2}\|p-q\|_1 \notag\\
  h^2(p,q) & \leq & \frac1{2} K(p,q).
\end{eqnarray}

The penalty to be defined shall depend on the size of the model, that we have
to estimate. The right tool is the metric entropy, and more precisely the
metric entropy with bracketing of the model.
\begin{dfntn}
\label{BE}
  Let $\mathcal{G}$ a function class. Let $N_{B,2}(\delta,\mathcal{G})$ be the smallest p such that there are couples of functions $[f_i^L,f_i^U]$ for i from 1 to $p$ that fulfill $\norm{f_i^L-f_i^U}_2\leq \delta$ for every j, and for any $f \in \mathcal{G}$, there is an $i \in [1,p]$ such that:
\begin{eqnarray*}
  f^L_i\leq f \leq f_i^U.
\end{eqnarray*} 
Then $H_{B,2}(\delta,\mathcal{G}) = \ln   N_{B,2}(\delta,\mathcal{G})$ is called the $\delta$-bracketing entropy of $\mathcal{G}$
\end{dfntn}
{\bf Remarks:} 
\begin{itemize}
\item{Notice that the $f_i^U$ and $f_i^L$ need not be in $\mathcal{G}$.}
\item{The $2$ in $H_{B,2}$ stands for $L^2$ distance.}
\end{itemize}

  Looking closely at definition \ref{BE}, we see that the concept of entropy
depends only on those of positivity and norms. We may then define a similar
bracketing entropy for any space with a norm and a partial order, such as the $L^1$ $\delta$-bracketing entropy of $\mathcal{Q}(N)$: we must find couples of Hermitian operators $[\tau_i^L,\tau_i^U]$ such that $\norm{\tau_i^U-\tau_i^L}_1\leq \delta$ and such that for any $\tau \in \mathcal{Q}(N)$, there is an i such that $\tau_i^L \leq \tau \leq \tau_i^U$.

We are chiefly interested in the $L^2$ entropy of square roots of density (denoted by $H_{B,2}(\delta,\mathcal{P}^{\frac1{2}})$):
\begin{eqnarray*}
\mathcal{P}^{1/2}(N) & = & \left\{ \sqrt{p_\rho} : p_\rho \in \mathcal{P}(N) \right\}. 
\end{eqnarray*}

Now the Theorem by Massart \cite{Mas}:

\begin{thrm}
\label{Massart}
Let $X_1,\dots,X_n$ be independent, identically distributed variables with
unknown density s with respect to some measure $\mu$. Let $(S_m)_{m \in \mathcal{M}}$
be an at most countable collection of models, where for each $m \in \mathcal{M}$, the
elements of $S_m$ are assumed to be densities with respect to $\mu$. We
consider the corresponding collection of maximum likelihood estimators
$\hat{s}_m$.
  Let $\pen : \mathcal{M} \longrightarrow \mathbb{R}$ and consider the random variable $\hat{m}$ such that:
\begin{eqnarray*}
 \prob{n}{-\ln(\hat{s}_{\hat{m}})} + \pen(\hat{m}) & = & \inf_{m\in \mathcal{M}}
\prob{n}{-\ln(\hat{s}_{m})} + \pen(m). 
\end{eqnarray*}

Let $(x_m)_{m \in \mathcal{M}}$ a collection of numbers such that
\begin{eqnarray*}
\sum_{m \in \mathcal{M}} e^{-x_m} & = & \sigma \leq \infty.
\end{eqnarray*}
For each m, we consider a function $\phi_m$ of $\mathbb{R}^{+*}$, nondecreasing, and such that  $x \mapsto \frac{\phi_m(x)}{x}$ is nonincreasing, and:
  \begin{eqnarray*}
 \phi_m(\sigma) & \geq & \int_0^{\sigma} \sqrt{H_{B,2}(\epsilon,S_m^{\frac1{2}})} d\epsilon. 
\end{eqnarray*}
We then define each $\sigma_m$ as the one positive solution of
\begin{eqnarray*}
 \phi_m(\sigma) & = & \sqrt{n} \sigma^2. 
\end{eqnarray*}
Then there are absolute constants $\kappa$ and C such that if for all $m \in \mathcal{M}$,
\begin{eqnarray*}
 \pen(m) & \geq & \kappa \left( \sigma_m^2 + \frac{x_m}{n}\right),  
\end{eqnarray*}
then
\begin{eqnarray*}
 \esp{h^2(s,\hat{s}_{\hat{m}})} & \leq & C \left( K(s,S_m) + \pen(m) + \frac{\sigma}{n} \right)
\end{eqnarray*}
where, for every $m \in \mathcal{M}$, $K(s,S_m) = \inf_{t \in S_m} K(s,t)$.
\end{thrm}

We notice that what is bounded \emph{in fine} is the Hellinger distance  and not the Kullback. Indeed our evaluation of the estimation error, which depends upon the size of the model (its bracketing entropy), dominates the Hellinger distance but maybe not the Kullback-Leibler distance.

In our case, we have parametrized the models $m$ by $N$, through definition \eqref{tronqueN}.

To apply Theorem \ref{Massart}, we need to find suitable $\phi_m$, and this calls for dominating the entropy integral. We reproduce here \cite{Gu.Gill.}.

By (\ref{comparehellinger}), it is sufficient to control $H_{B,1}(\delta,\mathcal{P}(N))$. Moreover, the linear extension of the morphism {\bf T} sends a positive matrix to a positive function, and is contractive. So any covering of $\mathcal{Q}(N)$ by $\delta$-brackets is sent upon a covering of $\mathcal{P}(N)$ by $L^1$ $\delta$-brackets, that is $[p^L_j,p^U_j]=[p_{\tau_j^L},
p_{\tau_j^U}]$. Thus 
\begin{eqnarray*}
H_{B,1}(\delta,\mathcal{P}(N)) & \leq & H_{B_1}(\delta,\mathcal{Q}(N)),
\end{eqnarray*}
so that
\begin{eqnarray*}
H_{B,2}(\delta,\mathcal{P}^{\frac1{2}}(N)) & \leq & C H_{B,1}(\delta^2,\mathcal{Q}(N)).  
\end{eqnarray*}
Moreover:
\begin{lmm}
\label{entropie}
\begin{eqnarray*} 
H_{B,1}(\delta,\mathcal{Q}(N)) & \leq & C N^2 \ln \frac{N}{\delta}
\end{eqnarray*}
where C is a constant not depending on $\delta$ or N, and can be put to  $1+\ln(5)$. 
\end{lmm}

\begin{proof}

Let $\left\{\rho_j: j=1,\dots,c(\delta,N) \right\}$ a maximal set of density matrices in $\mathcal{Q}(N)$ such that for all $j\neq k$, $\|\rho_j - \rho_k \|_1 \geq \frac{\delta}{2N}$. Define the brackets $[\rho_j^L,\rho_j^U]$ as 
\begin{eqnarray*}
\rho_j^L=\rho_j - \frac{\delta}{2N} {\bf 1} \quad \quad \rho_j^U = \rho_j +
\frac{\delta}{2N} {\bf 1}.
\end{eqnarray*}

Then $\|\rho_j^L - \rho_j^U\|_1 = \delta$. Moreover for any $\rho$ in the ball $B_1(\rho_j,\frac{\delta}{2N})$, as $\|\rho - \rho_j\|_1 \leq \frac{\delta}{2N} {\bf 1}$, we have
\begin{eqnarray*}
 \rho_j^L \leq \rho \leq \rho_j^U
\end{eqnarray*}
and as $\left\{\rho_j\right\}$ was a maximal set, this set of brackets cover $\mathcal{Q}(N)$. 

So $H_{B,1}(\delta,\mathcal{Q}(N)) \leq c(\delta,N)$.

Notice that $B_1(\rho_j,\frac{\delta}{4N})$ are disjoint and included in the shell  $B_1(0, 1 + \frac{\delta}{4N}) - B_1(0, 1 -\frac{\delta}{4N})$, so that

\begin{eqnarray}
\label{calcul}
 c(\delta,N) & \leq & \left(\frac{4N}{\delta}\right)^{N^2}\left( \left(1 +
\frac{\delta}{4N}\right)^{N^2} - \left(1 - \frac{\delta}{4N}\right)^{N^2} \right)
\notag \\
   & \leq & \left(1 + \frac{4N}{\delta}\right)^{N^2} \notag \\
   & \leq & \left( \frac{5N}{\delta}\right)^{N^2},
\end{eqnarray}
concluding the demonstration.

\end{proof}

From this, we can obtain:

\begin{crllr}
\label{entropiehellinger}
There is a constant $C$ such that:
$$H_{B,2}(\delta,\mathcal{P}^{\frac1{2}}(N)) \leq C N^2 \ln \frac{N}{\delta^2}.$$
 
\end{crllr}

Writing $$ \phi_N(\sigma)=\int_0^\sigma \sqrt{H_{B,2}(\epsilon,\mathcal{P}^{\frac1{2}}(N))} d\epsilon$$  and $\sigma_N(n)$ the only $\sigma$ such that
$$\phi_N(\sigma)=\sqrt{n}\sigma^2$$ we get
\begin{equation}
\label{sigN}
\sigma_N(n)\leq \sqrt{\frac{C}{n}}N \left( 1 + \sqrt{0 \vee \ln \frac{n}{N}}
\right).
\end{equation}
Indeed

\begin{eqnarray*}
   \phi_N(\sigma) & \leq & C N \int_0^\sigma \sqrt{\ln \left( \frac{N}{\epsilon^2} \right)} d\epsilon \\
               &  =   & C N^{\frac{3}{2}} \int^\infty_{\sqrt{\ln\frac{N}{\sigma^2}}} x e^{-\frac{x^2}{2}} dx \\
	       &  =   & C N^\frac{3}{2} \left( \int_{\sqrt{\ln\frac{N}{\sigma^2}}}^\infty e^{-\frac{x^2}{2}} dx 
                      - \left[ x e^{-\frac{x^2}{2}}\right]_{\sqrt{\ln\frac{N}{\sigma^2}}}^\infty \right) \\
  	       &  \leq & C N \sigma \left(1 + \sqrt{\ln \frac{N}{\sigma^2}} \right)
\end{eqnarray*} 
where we have made use of, in each  line in turn,
\begin{itemize}

\item{Corollary \ref{entropiehellinger}}
\item{the change of variables $x=\sqrt{\ln (N\epsilon^{-2})²}$, with $\frac{d\epsilon}{dx}=-\sqrt{N}x e^{-\frac{x^2}{2}}$}
\item{integration by parts, with $x$ seen as a primitive and $x e^{-\frac{x^2}{2}}$ as a derivative}
\item{the upper bound  $C e^{-\frac{x^2}{2}} $ of $\int_x^\infty
e^{-x^2/2} dx $ for $x$ positive when evaluating the first term.}
\end{itemize}

  We are looking for an upper bound on $\sigma_N$, solution of the equation

\begin{eqnarray*}
  \sqrt{n}\sigma_N^2 & = & C N \sigma \left( 1 + \sqrt{\ln \frac{N}{\sigma_N^2}} \right).
\end{eqnarray*}

We lower bound the second term by $0$, and get 
\begin{eqnarray*}
 \sigma_N & \geq & C \frac{N}{\sqrt{n}} \equiv \sigma_m 
\end{eqnarray*}

and upper bound
\begin{eqnarray*}
  \sigma_N & = & C N n^{-\frac1{2}}  \left( 1 + \sqrt{\ln \frac{N}{\sigma_N^2}} \right) \\
        & \leq & C N n^{-\frac1{2}}  \left( 1 + \sqrt{\ln \frac{N}{\sigma_m^2}} \right) \\
        & = & C \frac{N}{\sqrt{n}} \left( 1 + \sqrt{\ln \frac{n}{C^2 N}}
\right).   
\end{eqnarray*}

  We may absorb  the $C²$  in the first multiplicative constant to find (\ref{sigN}). Of course we take only the positive part of the logarithm. This will always be the case hereafter.

Applying Theorem \ref{Massart} we get:

\begin{thrm}

Consider the collection of maximum likelihood estimators $(\hat{\rho}_N)_{N\in \mathbb{N}}$, that is for any integer N,
\begin{eqnarray*}
 \prob{n}{- \ln(p_{\hat{\rho}_N)}} & = & \inf_{\rho \in \mathcal{Q}(N)} \prob{n}{- \ln(p_{\hat{\rho})}}
\end{eqnarray*}
Let $ \pen: \mathbb{N} \mapsto \mathbb{R}_+$ and consider a random variable  $\hat{N}$ such that 
\begin{eqnarray*}
\prob{n}{- \ln(p_{\hat{\rho}_{\hat{N}}})} + \pen(\hat{N}) & = & \inf_{N \in \mathbb{N}} 
(\prob{n}{- \ln(p_{\hat{\rho}_N})} + \pen(N))
\end{eqnarray*}
Let $(x_N)_{N\in \mathbb{N}}$ a family of positive numbers such that 
\begin{eqnarray*}
 \sum_{N\in \mathbb{N}} e^{-x_N} & = & \sigma  ~<~ \infty
\end{eqnarray*}
Then there are absolute constants  $\kappa$ and C such that if 
\begin{eqnarray*}
\pen(N) & \geq & \kappa
(\frac{N^2}{n}(1 + (0 \vee \ln\frac{n}{N})) + \frac{x_N}{n})
\end{eqnarray*}
then 
\begin{eqnarray*}
 \mathbb{E}[h^2(p_{\rho}, p_{\hat{\rho}_{\hat{N}}})] & \leq & C \left( \inf_{N \in \mathbb{N}} (\mathbb{E}[K(\rho, \mathcal{Q}(N))] + \pen(N)) + \frac{\Sigma}{n} \right) 
\end{eqnarray*}
with $K(\rho, \mathcal{Q}(N)) = \inf_{\tau \in \mathcal{Q}(N)} K(p_{\rho}, p_{\tau})$.
 \end{thrm}

\smallskip
{\bf Remarks:} 
\begin{itemize}
{\item
 When designing the penalty, what stands out in this theorem is the general
form of the penalty. Now the constant $\kappa$ that can be explicitly computed
would be very pessimistic. The best thing to do is therefore to keep the
general formula for the penalty and calibrate $\kappa$ using cross-validation,
the slope heuristic \cite{Mas} or any other appropriate method. 
}
{\item
  If we wanted an explicit convergence rate for a given state, as for the photon
basis in section \ref{photon}, we would first need to know how the
Kullback-Leibler distance $K(\rho,\mathcal{Q}(N))$ is decreasing with $N$. One thing
that is obvious, however, is that if we add noise we convolve with the same
function $p_{\rho}$ and $p_{\sigma}$ for all $\sigma$ in $\mathcal{Q}(N)$, so the
Kullback-Leibler distance is decreasing with the noise, so convergence is
faster when there is noise... The reason for this is that we are looking at
convergence in Hellinger distance, that is a distance between the law of the
result of the measurement $p_{\rho}$ and $p_{\sigma}$. This does not tell us
directly anything about what we are really interested in, that is the distance between
$\rho$ and $\sigma$ (as operators). Indeed we may bound the $L^2$ or $L^1$
norm between elements of $\mathcal{Q}(N)$ by the Hellinger distance, times something
depending on the sum of the $L^2$ or $L^{\infty}$ norms of the
$f_{j,k}^{\eta}$. And these norms are going (very fast) to infinity when there
is noise, so that low Hellinger distance gives no indication on the operator
norms. 
}
\end{itemize}

\section{Quantum calibration of a photocounter}
\label{photocounter}

 This section features a scheme to calibrate an apparatus $M$ measuring the
number of photons in a beam with the help of a photocounter. 

The physical motivation is given in Appendix \ref{phy.photo}.

The first subsection states the mathematical problem.
In the two others are studied
respectively projection estimators and maximum likelihood estimators.

\subsection{Statistical problem}

The practical problem of calibration of a photocounter turns out to be mathematically speaking an entirely classical missing data problem. However, to the best of our knowledge, it has never been studied. We now describe this missing data problem.

We are given samples $(i,x)$ in $\mathbb{N} \times \mathbb{R}$ from a probabiliy density of the form 
\begin{eqnarray}
\label{lawp}
p(i,x) & = & \sum_{k=0}^{\infty} b_k^2 P^i_k \psi_k(x)^2. 
\end{eqnarray}
In this expression, the real numbers $b_k^2$ satisfy $\sum_m b_k^2 = 1$. The $\psi_k$ are the Fock basis functions given in Equation \eqref{hermite}. For any $k$, the $P^k_i$ are a probability measure, that is they are non-negative and $\sum_{i=0}^{\infty} P_i^k = 1$.

We know the $b_k^2$, and we want to retrieve the $P_i^k$, which we do not know. We write $P = (P_i^k)_{i,k}$.

To make clearer that this is a missing data problem, we give the following way to obtain this experiment.
First we choose $k\in \mathbb{N}$  with probability given by $b_k^2$. We forget $k$, which is the missing data. Our data consists in  $(i,x)$, with $i$ having law given by  $P_i^k$ and $x$ with law $\psi_k(x)$.

Notice that the experimentalist has some control on the $b_k^2$, but usual techniques will yield $b_k^2$ proportional to $\xi^k$. This means that the low $k$ are probed faster.

\smallskip

We propose below two types of estimators $\hat{P}$ for $P$. To get results on their efficiency, we must first find  meaningful distance $d(P,\hat{P})$. 
Since $\sum_i P_i^k =1$ for all $k\in \mathbb{N}$,
distances like $d_2^2(P,Q) = \sum_{i,k} (P_i^k - Q_i^k)^2$ are bound to yield
infinite errors on our estimators. We then must weight them, using $(a_k)_{k\in \mathbb{N}}$ of our choice. We shall use,
depending on the estimator, either $d_2^2(P,Q) = \sum_{i,k} a_k^2 (P_i^k -
Q_i^k)^2$ with $\sum a_k^2 = 1$, or $d_1(P,Q) = \sum_{i,k} a_k |P_i^k -
Q_i^k|$, with $\sum_{k} a_k = 1$. Then these distances are bounded by $2$ on
the set of all $P$ such that  $\{P_i^k\}_{i\in \mathbb{N}}$ is  a probability measure for every $k$. 

 Varying the choice of $a_k$ corresponds to putting the emphasis on different $k$, that is deciding which $P_i^k$ we demand to know with the more precision. If we  take the $a_k$ decreasing, it means physically that we are more interested in the
behaviour of our photocounter for a low number of photons. This is usually
the case for a physicist. A possible choice is to take $a_k$ or $a_k^2$ equal to $b_k^2$.

In the next subsection, we use projection estimators, and in the following, maximum likelihood estimators.

\subsection{Using projection estimators}

 As in the tomography problem,  the parameter space is contained in an infinite-dimensional vector
space, and a natural type of estimators are projections of the empirical law on
finite-dimensional subspaces. The problem we are left with is then again
finding the best subspace.

Concretely, we consider the distance $d_2^2(P,Q) = \sum_{i,k} a_k^2 (P_i^k -
Q_i^k)^2$ and write $E_i^k = a_k P_i^k$. Similarly we shall write $\hat{E}_i^k = a_k \hat{P}_i^k$ for our estimator. Then
\begin{eqnarray*}
d^2_2(P, \hat{P}) & = & \sum_{i,k} (E_i^k - \hat{E}_i^k)^2,
\end{eqnarray*} 
and the law of our samples can be rewritten as 
\begin{eqnarray}
\label{law1}
p(i,x) & = & \sum_k E_i^k \frac{b_k^2}{a_k} \psi_k(x)^2.
\end{eqnarray}
We may then consider $\{(b_k^2/a_k) \psi_k \mathbf{1}_{i=l}\}_{k,i}$ as a basis of our functions on $\mathbb N \times \mathbb R$. We want to use the general constructions of section \ref{projection}. We first need a dual basis $\{g_{i,k}\}$. Now, the dual basis of $\{\psi_k^2\}$ as functions on $\mathbb R$ is well-known. Those are the ``pattern functions'' $f_{k,k}$ introduced in \cite{D'Ariano.0} (see \eqref{pattern}). From this, we deduce:
\begin{eqnarray*}
g_{i,k}(l,x) & = & \frac{a_k}{b_k^2} f_{k,k}(x) \mathbf{1}_{i=l}.
\end{eqnarray*}
With these dual functions, we can define the minimum contrast function:
\begin{eqnarray*}
\gamma_n( Q) & = & d_2^2(Q, 0) - 2 \left(\sum_{\alpha=1}^n \frac{g_{i,k}(l_{\alpha}, x_{\alpha})}{a_k}\right) \left(\sum_{i,k} a_k^2 Q_i^k\right)
,
\end{eqnarray*}
where the $(l_{\alpha}, x_{\alpha})$ are our data, that is $n$ independent samples with law $p$.

Our models $m\in\mathcal M$ consist in the subsets of $\mathbb{N}^2$. If $(i,k)\not\in m$, then $\hat{P}^k_i = 0$. In a model $m$, the estimator $\hat{P}^{(m)}$ given by minimizing the contrast function is then 
\begin{eqnarray*}
\hat{P}_i^k & = & \frac1n \sum_{\alpha=1}^{n} \frac{g_{i,k}(l_{\alpha},x_{\alpha})}{a_k} \mathrm{ ~~for~~ } (i,k)\in m. 
\end{eqnarray*}

The penalized estimator is as always the projection estimator of the model $\hat{m}$ such that:
\begin{eqnarray*}
\hat{m} & = & \arg\min_{m \in \mathcal{M}} \gamma_n( \hat{P}^{(m)} ) +
\pen_n(m).  
\end{eqnarray*}
We also use the usual notation for the distance to a model:
\begin{eqnarray*}
d_2(P,m) & = & \inf_{Q \in m} d_2(P,Q).
\end{eqnarray*}  

We then obtain from the general theorems of section \ref{projection}:

\begin{thrm}
\label{phhoeff}
  Let $P$ be a photocounter and $(a_k)$ and $(b_k)$ with $\sum_k a_k^2 = \sum_k b_k^2 = 1$. Let $(x_{i,k})_{(i,k)\in \mathbb{N}^2 }$ such that $\sum_{i,k} e^{-x_{i,k}} = \Sigma<\infty$. We define a penalty as 
\begin{eqnarray*}
\pen_n(m) & = & \sum_{(i,k)\in m} (1+ \epsilon) \left( \ln(M_{i,k}) + \frac{x_{i,k}}{2}  \right) \frac{M_{i,k}^2}{n} 
\end{eqnarray*}
with
\begin{eqnarray*}
M_{i,k} & = & \frac{a_k}{b_k^2}(\sup_x f_{k,k}(x) - \inf_x f_{k,k}(x)).
\end{eqnarray*}

Then the penalized estimator fulfills 
\begin{eqnarray*}
\mathbb{E}\left[\frac{ \epsilon }{2 + \epsilon }  d_2^2(P, \hat{P} )\right]
&  \leq & \inf_{m\in \mathcal{M} } \left(1+ \frac{2}{ \epsilon
}\right)d_2^2(P,m) + 2 \pen_n(m) + \frac{(1+ \epsilon) \Sigma  }{n}.
\end{eqnarray*} 
\end{thrm} 

\begin{thrm}
\label{phbern}
 Let $P$ be a photocounter and $(a_k)$ and $(b_k)$  with $\sum_k a_k^2 = \sum_k b_k^2 = 1$. Let $(y_{i,k})_{(i,k) \in \mathbb{N}^2 }$ such that $\sum_{i,m} e^{-y_{i,m}} = \Sigma<\infty$. Let then
\begin{eqnarray*}
x_{i,k} & = & 2 \ln\left(\frac{a_k}{b_k^2} \left\lVert f_{k,k}
\right\rVert_{\infty} \right) + y_{i,k}.
\end{eqnarray*} 

  For any $\delta \in (0,1)$, with 
\begin{eqnarray*}
\pen_n(m) & = & \sum_{(i,k)\in m} \pen_n^{(i,k)} \\
\pen_n^{(i,k)}&  = & 2 \frac{1+ \epsilon }{n} \left( \sqrt{\frac{2}{1- \delta }
x_{i,k} \left( \mathbb{P}_n[g_{i,k}^2] + \frac 1n
\frac{a_k^2}{b_k^4}\left\lVert f_{k,k}\right\rVert^2_{\infty} \left(\frac 13 +
\frac 1 \delta \right)  x_{i,k}  \right)  } + \frac{ a_k \left\lVert f_{k,k}
\right\rVert_{\infty} }{3b_k^2\sqrt{n}} x_{i,k} \right)^2, 
\end{eqnarray*}
there is a constant $C$ such that:  
\begin{eqnarray*}
\mathbb{E}\left[\left(\frac{ \epsilon }{2 + \epsilon } d_2^2(P, \hat{P})  - \left( \left(1 + \frac 2 \epsilon \right) \inf_{m \in \mathcal{M}_n} d^2_2(P,m) + 2 \pen_n(m) \right)\right) \vee 0 \right] & \leq & \frac{C \Sigma }{n}
\end{eqnarray*}
where $ \mathcal{M}_n$ is the set of models $m$ for which $(i,k) \in m $ implies $ x_{i,k} < n$.

\end{thrm} 

\smallskip
{\bf Remarks:} 
\begin{itemize}
\item{
As with the estimation of states with tomography in section \ref{projection}, we choose with high efficiency the best subspace. It should be noticed
that convergence is fast if the photocounter is good, and could be slower
if it is bad. In the latter case, we know it is bad, though. Indeed, the
dependence of the convergence rate on the photocounter $P$  lies in the
approximation properties of the models -- subspaces -- $m$, that is on how fast  $d^2_2(P,m)$ decrease when $m$ gets bigger. Now for an ideal
photocounter, we need only the $(i,i)$ to be in $m$. The penalty would be as
low as possible when neglecting what happens to beams with more than a
given number $k$ of photons. For a worse photocounter,
to have a good approximation of how a $k$-photons beam is read,
we might need many $i$, and the penalty would include all the $\pen^{i,k}$.
}
\item{
The estimator depends only weakly on $(a_k)$ (unlike the distance), which is
good news as it is somewhat arbitrary. Indeed, the empirical $ \hat{P}^k_i$
does not depend of this sequence at all, nor do the main terms in the threshold
on $ \hat{P}^k_i$  of both theorems. For Theorem \ref{phhoeff}, this main
term is $\sqrt{a_k^{-1}(1 + \epsilon ) \ln(B_{i,k})} B_{i,k}/\!\sqrt{n}$. Now
$B_{i,k}$ depends linearly on $a_k$, so the only $a_k$ left in this expression
is in the logarithm which can be developed as $\ln(B_{i,k}/\!a_k) + \ln(a_k)$.
In this way, we see that we only get another term in the penalty. For Theorem \ref{phbern}, the threshold is essentially $a_k^{-1} \sqrt{8(1+ \epsilon ) \mathbb{P}_{n}\left[g_{i,k}^2\right]\ln( \left\lVert g_{i,k}\right\rVert_{\infty} ) /\!((1- \delta)n )}$; and as $g_{i,k}$ is proportional to $a_k$, the situation is the same.   
}
\item{ The process by which we get our data includes a tomographer and the laws $p(i,x)$ were given in the ideal case when there is no noise. If there is noise, as briefly sketched in section \ref{noise}, these laws are different. However we may characterize the noise with a single $0 < \eta < 1$.
We then have for free the same theorems for $ \eta > \frac 12$: we only need to replace $f_{k,k}$ with $f^{\eta}_{k,k}$.
} 
\end{itemize}

\subsection{Maximum likelihood procedure} 

  In this case,  our results  are easier expressed with the distance 
\begin{eqnarray*}
 d_1(P, \hat{P}) & = & \sum_{i,k} a_k \left| P_i^m - \hat{ P}^k_i \right| \\
                     & = & \sum_{i,k}     \left| E_i^k - \hat{ E}^k_i \right|
\end{eqnarray*} 
with $E_i^k = a_k M_i^k$ and $\sum_k a_k =1$. We  denote $w_i = \sum_k E^k_i$. Notice that $\sum_i w_i = 1$. 

Recall that our data consists in $n$ independent samples $(l_{\alpha}, x_{\alpha})$ with law $p$ given by Eq. \eqref{lawp}.

  The main difficulty with applying here Theorem \ref{Massart}  lies in that  the
Kullback distance to the models is usually infinite (if we have
$\hat{E}_i^k=0$ for all $k$ for some $i$, then $\hat{p}(i,\mathbb{R})=0$ and
this is generally not the case for $p(i,\mathbb{R})$). The easiest way around
is to keep independence and restrict attention to some set of $i$.

  Explicitly, we take an ordering on the possible results $i$ of the photocounter
(typically, if we expect that one result  corresponds roughly to a given
number of photons, we can order them in increasing order. The idea is that
the results that interest us most should come first). 
We then choose, still  beforehand,  $ I_{max}\in\mathbb{N}$, and we
restrict our attention  to the first $  i\in [0, I_{max}] $. We just throw away the part of the data where the photocounter
gave a result more than $ I_{max}$. We are left with data size $n_{
 I_{max} }$, with law $p_{  I_{max} }$ on $  [0,I_{max}] \times\mathbb{\mathbb{R}}$:
\begin{eqnarray*}
p_{  I_{max} } & = & \frac{p_{| [0,I_{max}] \times\mathbb{R}}}{\int_{
 [0,I_{max}] \times\mathbb{R}}p} .
\end{eqnarray*}

This law is the probability measure associated to the apparatus $\tilde{P}$ for which
$\tilde{P}_i^k = \frac{1}{\sum_{l \leq  I_{max}  }w_l} P^k_i \mathbf{1}_{i \leq
 I_{max}  }$.

  The models $m_{ I ,K}$ we work with are indexed by
$K\in
\mathbb{N}  $ and $  I  \leq I_{max}  $. They are given by the constraints:
\begin{align}
  \hat{E}_i^k & = 0 & \mbox{if }& i >  I_{max}  \notag \\
\hat{E}_i^k & = 0 & \mbox{if }& i >  I_{max}  \mbox{ and } k \leq K
 \notag \\
  \sum_{i \leq  I  } \hat{E}_i^k &  = a_k & \mbox{for } & k \leq K
  \notag \\
\label{condition4} 
 \hat{E}_i^k & = \frac{a_k}{ I_{max} + 1 } & \mbox{for } & k > K  \mbox{
and } i \leq I_{max}.  
\end{align} 

  Any such element gives a probability measure on $ ([0, I_{max}] \times \mathbb{R})$.
Similarly to equation (\ref{law1}), the corresponding probability law reads $ \hat{p}(l,x) = \sum_{i,k} b_k^2 a_k^{-1}
\hat{E}_i^k \psi_k(x)^2 \mathbf{1}_{i=l}$. The fourth condition
(\ref{condition4}) does not increase the complexity of the model and ensures
that the Kullback distance remains finite.

We can now use an empirical maximum likelihood procedure to select within each model an estimator. It  minimizes on each $m_{I,K}$ the contrast function
\begin{eqnarray*}
\gamma_n( Q) & = & \sum_{k=1}^n - \ln q(l_k,x_k). 
\end{eqnarray*}
where $Q$ is an element of the model $m_{I,K}$ and $q$ the associated probability law.

  We then use Theorem \ref{Massart} to select the model of which we keep the estimator, through a penalization procedure. We obtain the following theorem.
\begin{thrm}
  Consider the collection of maximum likelihood estimators $
(\hat{P}_{  I,K })_{I \leq
 I_{max},  K \in \mathbb{N}}$, defined as minimizers of 
\begin{eqnarray*}
\gamma_n( \hat{P}_{  I , K }) & = & \inf_{P \in m_{
 I ,K }} \gamma_n(P)  
\end{eqnarray*}
Let $\pen: [0, I_{max} ]\times\mathbb{N} \rightarrow \mathbb{R}$ be a
penalty function and define $ (\hat{  I  }, \hat{K})$ by 
\begin{eqnarray*}
\gamma_n( \hat{M}_{ (\hat{  I  }, \hat{ K})}) + \pen(
\hat{  I  }, \hat{ K}) & = & \inf_{  I \leq
 I_{max}  , K  \in \mathbb{N} }  \gamma_n( \hat{P}_{
 I , K }) + \pen(  I , K ). 
\end{eqnarray*}
Let $(x_{  I , K })$ be a family of numbers such that 
\begin{eqnarray*}
\sum_{  I \leq  I_{max}  , K  \in \mathbb{N} } e^{- x_{
 I , K }} & = & \Sigma ~<~ \infty.
\end{eqnarray*}
  Then there are absolute constants $ \kappa $ and $C$ such that if
\begin{eqnarray*}
\pen(  I  , K) & \geq & \kappa \left(  (I+1)(K+1)
\frac{\ln(n_{  I_{max} })}{n_{  I_{max}  }} + \frac{x_{  I  ,
K }}{n_{  I_{max}  }}\right), 
\end{eqnarray*}
then
\begin{eqnarray*}
\mathbb{E}\left[ d_1(P, \hat{P}_{ (\hat{  I  }, \hat{K})
}) \right] & \leq & \sum_{i > I_{max}  } w_i + \sum_{k \in \mathbb{N}} \left(
2 a_k \wedge \left(C  \frac{a_k}{b_k^2} \left\lVert f_{k,k}
\right\rVert_{\infty} \sqrt{\inf_{  I \leq  I_{max}   ,K
\in \mathbb{N} } K(p_{  I_{max} }, m_{  I , K })+\pen(
 I , K ) + \frac{\Sigma}{n_{ I_{max} }} }  \right) \right),
\end{eqnarray*}
where $K(p_{  I_{max} }, m_{  I , K }) = \inf_{Q \in m_{
 I , K}} K(p_{  I_{max} }, q)$, intended as the Kullback
distance  on $ [0,I_{max}] \times \mathbb{R}$.

\end{thrm}

{\bf Remarks:}
\begin{itemize}
{\item
As with projection estimators, we can expect fairly quick approximation if the
photocounter is good. Indeed, for $K = I_{max} $ and the ideal photocounter, the distance
$K(p_{ I_{max} },m_{  I_{max}, K })= 0$.
}
{\item
 Like projection estimators, the maximum likelihood strategy can also be used with noise. If
$\eta > \frac 12$, we get the same theorem changing $f_{k,k}$ in
$f_{k,k}^{\eta}$. Just notice that the infinite norm $\|f_{k,k}\|_{\infty}$ is exploding.
}
{\item
As in section \ref{EMV}, an explicit computation of $\kappa$ would be
over-pessimistic and it is best to estimate it with a data-driven procedure.
}
\end{itemize}

\begin{proof}
  First we rewrite and bound the distance $d_1$ in a way that suits our
purpose. We separate the entries corresponding to measurement results bigger
than  $
 I_{max} $, and we recall at the third line that $\sum_{i \in \mathbb{N} } E^k_i = a_k$.  Then 
\begin{align*}
d_1(P, \hat{P}) & = \sum_{i,k} \left| E_i^k - \hat{E}_i^k   \right| \\
                    & = \sum_{i > { I_{max} } } \sum_k E^k_i + \sum_k
\sum_{i \leq { I_{max} }}  \left| \hat{E}_i^k - E^k_i \right|   \\
		    & \leq \sum_{i > I_{max} } \sum_k E^k_i + \sum_k
\left( 2 a_k\wedge \left(\sum_{i \leq { I_{max} }}  \left| \hat{E}_i^k -
\frac{1}{\sum_{i \leq  I_{max}}  w_i} E^k_i     \right| + \left(\frac{1}{\sum_{i
\leq  I_{max}  } w_i} - 1 \right) E^k_i  \right) \right) \\  
                    & =   \sum_{i >  I_{max}  } w_i + \sum_{i \leq
I_{max} } \frac{\sum_{i >I_{max} } w_i}{\sum_{i \leq  I_{max} }w_i
}\sum_k E_i^k + \sum_k \left(2 a_k \wedge \sum_{i \leq I_{max}} \left|
\hat{E}_i^k - \frac{1}{\sum_{i \leq I_{max}}  w_i} E^k_i     \right| \right) \\
                    & =  2  \sum_{i >  I_{max}  } w_i + \sum_k \left( 2
a_k \wedge \sum_{i \leq I_{max}} \left| \hat{E}_i^k - \frac{1}{\sum_{i \leq
 I_{max} } w_i} E^k_i     \right| \right). \\ 
\end{align*}

Let us now work a little on the last term:
\begin{eqnarray*}
\frac{1}{\sum_{i \leq  I_{max} } w_i} E_i^k  & = & \int \frac{a_k}{b_k^2}
f_{k,k}(x) \mathbf{1}_{i=l} dp_{  I_{max}  }(l,x), \\
\hat{E}^k_i                               & = & \int \frac{a_k}{b_k^2} f_{k,k}(x) \mathbf{1}_{i=l} d\hat{p}(l,x).  
\end{eqnarray*}
So that 
\begin{eqnarray*}
\left| \frac{1}{\sum_{i \leq I_{max} } w_i } E_i^k  - \hat{E}^k_i\right| &
= & \left| \int f_{k,k}(x) \mathbf{1}_{i=l} d(p_{  I_{max} } - \hat{p})(l,x)  \right| \\
									& \leq &
\frac{a_k}{b_k^2} \left\lVert f_{k,k} \right\rVert_{\infty} \int
\mathbf{1}_{i=l} d|p_{  I_{max} } - \hat{p}|(l,x).   
\end{eqnarray*}

Summing over $i$, we get:
\begin{eqnarray*}
\sum_{i \leq  I_{max}  } \left| \frac{1}{\sum_{i \in  I_{max}} w_i }
E_i^k  - \hat{E}^k_i\right| 	& \leq & \frac{a_k}{b_k^2}  \left\lVert f_{k,k} \right\rVert_{\infty} \int
d|p_{  I_{max} } - \hat{p}|(l,x).  
\end{eqnarray*} 

We may then bound the distance between the POVM we calibrate and our estimator by
\begin{eqnarray*}
d_1(P, \hat{P}) &=& 2 \sum_{i >  I_{max} } w_i + \sum_{k \in \mathbb{N}}
\left( 2 a_k \wedge \left( \frac{a_k}{b_k^2} \left\lVert f_{k,k}
\right\rVert_{\infty}  \int
d|p_{  I_{max} } - \hat{p}|(l,x)  \right) \right)   .
\end{eqnarray*}

  Finishing the proof of our theorem amounts to controlling $ \int
d|p_{  I_{max} } - \hat{p}|(l,x)
$. We first apply Theorem \ref{Massart} (assuming that our penalty is big enough, which we check below).   We get:
\begin{eqnarray*}
\mathbb{E}\left[h^2(p_{  I_{max} } , \hat{p}_{(\hat{  I  },
\hat{ K} )} )  \right] & \leq & C \left(\inf_{  I \leq
 I_{max} ,K  \in
\mathbb{N} } K(p_{  I_{max}}, m_{  I , K })+\pen(
 I,K) + \frac{\Sigma}{n_{ I_{max} }}\right).
 \end{eqnarray*} 
We then use the bound (\ref{comparehellinger})  of the square of the $L^1$-distance in the Hellinger distance, and finish with Jensen, using the concavity of both the function $x \mapsto (C \wedge x)$ and the square root.
\begin{eqnarray*}
\mathbb{E}\left[ d_1(P, \hat{P}_{ (\hat{  I  }, \hat{K})
}) \right] & \leq & \mathbb{E}\left[ \sum_{i > I_{max}  } w_i + \sum_{ k
\in \mathbb{N}} \left( 2 a_k \wedge \left(C  \frac{a_k}{b_k^2} \left\lVert
f_{k,k}  \right\rVert_{\infty} \int
d|p_{  I_{max} } - \hat{p}_{ (\hat{  I  }, \hat{K})}|(l,x)
\right)\right) \right] \\
 & \leq & \sum_{i > I_{max}} w_i + \sum_{k \in \mathbb{N}}
\mathbb{E}\left[  \left( 2 a_k \wedge \left(C  \frac{a_k}{b_k^2} \left\lVert
f_{k,k}  \right\rVert_{\infty}  \sqrt{h^2\left( p_{  I_{max} } -
\hat{p}_{ \hat{  I  }, \hat{ K}}  \right)}  \right)\right)  \right] \\
 & \leq & \sum_{i > I_{max}  } w_i + \sum_{k \in \mathbb{N}} \left( 2 a_k
\wedge \left( C \frac{a_k}{b_k^2} \left\lVert f_{k,k}  \right\rVert_{\infty}
\sqrt{ \mathbb{E}\left[ h^2\left( p_{  I_{max} } - \hat{p}_{ \hat{
 I  }, \hat{ K }}  \right)  \right]}\right)\right) \\
 & \leq & \sum_{i > I_{max}  } w_i + \sum_{k \in \mathbb{N}} \left( 2 a_k
\wedge \left(C  \frac{a_k}{b_k^2} \left\lVert f_{k,k}  \right\rVert_{\infty}
\sqrt{\inf_{  I \leq  I_{max} , K \in \mathbb{N} } K(p_{  I_{max} }, m_{  I ,
K })+\pen(  I , K ) + \frac{\Sigma}{n_{ I_{max} }} }  \right) \right). 
\end{eqnarray*}

The only thing we still have to check is our penalty. We must dominate $
H_{B,2}( \delta, \mathcal{P}^{1/\!2}(  I , \mathcal{M} ))$ where 
\begin{eqnarray*}
 \mathcal{P}^{1/\!2}(  I , K ) & = & \left\{ \sqrt{q}, Q \in
m_{  I , K } \right\}. \\
\end{eqnarray*} 

 With the same reasoning as in section \ref{EMV}, it is sufficient to dominate
$H_{B,1}( \delta^2, m_{  I , K })$. We then mimic lemma
\ref{entropie}. All the elements of $m_{  I , K }$ are on the
$L^1$-sphere of radius $\sum_{k \leq K }a_k$ of a vector space of
dimension $ (K + 1)   (I + 1) $. We can then associate a maximal
collection of brackets to a maximal collection $(P_j)$ of $P \in m_{
 I , K }$ separated by $ \delta^2 /\!(2 (K+1)
 (I+1)   )$. The balls $B_1(M_j, \frac{ \delta^2}{ (K+1)
 (I+1) })$ are disjoint and in the shell $ B_1(0, \sum_{k\leq K }
a_k +  \frac{ \delta^2}{ (K+1)    (I+1) }) - B_1(0, \sum_{k \leq K }a_k - \frac{ \delta^2}{(K+1) (I+1) }) $.   And as with equation (\ref{calcul}), we obtain 

\begin{equation*}
  H_{B,1}(\delta^2, m_{  I , K }) \leq C (K+1)    (I+1)  \ln\left(\frac{
(K+1)(I+1) }{\delta^2} \right) 
\end{equation*} 

  Imitating the calculation in the proof of corollary \ref{entropiehellinger},
we find that the solution $ \sigma_{  I , K }$ of the equation
\begin{eqnarray*}
  \sqrt{n_{  I_{max} }}\sigma_{  I , K }^2 = \int_0^{ \sigma_{
 I , K }} \sqrt{H_{B,2}(\delta, \mathcal{P}^{1/\!2}(
 I , K )})
\end{eqnarray*}
admits this upper bound:
\begin{eqnarray*}
\sigma_{  I , K } \leq C \sqrt{\frac{ (K+1)
 (I+1) }{n_{  I_{max} }}} (1 + \sqrt{\ln n_{  I_{max} }})
\end{eqnarray*}
We may absorb the latter $1$ in the constant, as long as $n_{I_{max}} \geq 2$...

This ends the proof. 

\end{proof}


\section*{Acknowledgements}

The main part of this work stems from my master's thesis, written under the
very pleasant direction of Pascal Massart. Patricia Reynaud was also most
helpful during this period.

I am equivalently indebted to Richard Gill and to Madalin Gu\c{t}\u{a} in
Eindhoven for useful discussions.

The presentation of this paper has been greatly improved and several errors in
the proofs suppressed (thanks Richard!) by their careful rereading, together
with that of Cristina Butucea.


\appendix

\section{Background in quantum mechanics}
\label{Quantique}

Subsection \ref{Formalism} gives parallel developments of classical statistics and
quantum statistics, so that any quantum notion is linked with a classical
equivalent. 

Subsection \ref{QHD} describes both the experimental setup of
quantum homodyne  tomography and some basic mathematics playing a role in it.
More precisely, it  highlights several different representations of the state
to be recovered (our  unknown) and the links between them. 

Subsection \ref{phy.photo} is background for section \ref{photocounter}. Notably, it explains where the formulas such as \eqref{law1}  come from.

\subsection{Statistics: classical and quantum}
\label{Formalism}

We have here three different parts. The aim is to highlight the equivalences
in classical and quantum formalism. The first part lies then upon the
classical world, the second part recast this construction as a special
case of what will be our quantum formalism, and the third part describes these
quantum statistics. Bold numbers refer to the same number in the other
sections.
They might be repeated inside a section if the same object is introduced
under different forms.

In this short introduction to the subject, we shall restrict ourselves more or less to describing what physical
measurements can be done and how they can be encoded mathematically. In other
words, we characterize what information can be retrieved from a system.

\subsubsection{Classical}
\label{classique}

In the classical setting of statistics, we are working with probability measures $p$\
\rv{1}
on a probability space $(\mathcal{X},\mathcal{A})$\ \rv{2}. For comparison, we recall that probability measures
are normalized \rv{3} real \rv{4} non-negative \rv{5} measures. Similarly measures
are  elements of
$\mathcal{M}(\mathcal{X},\mathcal{A})$\ \rv{6}, the dual of $L^{\infty}(\mathcal{X},\mathcal{A})$\ \rv{7}.

Notice that the probability measures form a convex set, the extremal points of which
are the Dirac measures \rv{8} on $x$ for $x\in (\mathcal{X},\mathcal{A})$. They may then be
described by $x$\ \rv{9}. If we want to draw on the analogy with physics
$(\mathcal{X},\mathcal{A})$ may be viewed as a phase space, and the $x$ would be the pure
states. A general probability measure would describe a mixed state. These are systems that have
a probability to be in this or that pure state. Any mixed state (probability measure)
can be decomposed in a unique way over pure states (Dirac).

A statistical \emph{model}\ \rv{10} consists in a set of probability measures
$p_{\theta}$ on a probability space $(\mathcal{X},\mathcal{A})$ indexed by a  parameter
$\theta$, for $\theta\in \Theta$\ \rv{11} the \emph{parameter space}. A
statistical problem consists in determining as precisely as possible, with a
meaning depending on the instance, a function of $\theta$. 

Now we must gain access at information on  these $\theta$ in some way. What we
have access at are random variables. 

The aforementioned space $L^{\infty}(\mathcal{X},\mathcal{A})$ is the space of real bounded random
variables $f$\
\rv{12}. By analogy with the quantum case, we call these
$f$ \emph{observables}. They correspond to the set of physical measurements
that can be carried out on the system, to what can be ``observed''.

``Measuring'' an observable $f$  yields a result $f(x)$\ \rv{13}, with law:
\begin{align}
\label{lawresult}
\prob{p}{f \in B} & = \int_{\mathcal{X}} {\bf 1}_{f(x)\in B} dp(x) & \mbox{for\ } B \in
\mathcal{B} \mbox{\ $\rv{14}$}
\end{align}
where $\mathcal{B}$ is the borelian $\sigma$-algebra of $\mathbb{R}$.
Notice that this result is not random for a pure state.

Notice also that the way we could see the probability measures $p$ as elements of the
dual of
$L^{\infty}(\mathcal{X},\mathcal{A})$ was by writing  $p(f) = \int_{\mathcal{X}} f(x)\dd p(x)$\
\rv{15}.

The most general type of statistic or estimator we can extract from data,
including random strategies, is obtained by associating to each $x$ a
probability measure on an auxiliary space $(\mathcal{X}_a, \mathcal{A},a)$\ \rv{16} and draw a final result
according to this probability measure. This is equivalent (at the price of changing the
auxiliary space) to measuring a function $f$\ \rv{17} on a space
$(\mathcal{X}\otimes\mathcal{X}_a,\mathcal{A}\otimes\mathcal{A}_a)$\ \rv{18} according to a probability measure
$p_{\theta}\otimes s$\ \rv{19} with $s$ independent of $\theta$.

If we write (\ref{lawresult}) in this case, we get
\begin{align*}
\prob{\theta}{f \in B} & =  \int_{\mathcal{X}}\int_{\mathcal{X}_a} {\bf 1}_{f(x,x_a)\in B}
dp_{\theta}(x) ds(x_a) & \mbox{for\ } B\in
\mathcal{B}. 
\end{align*}

If we integrate out $\mathcal{X}_a$, this yields
\begin{align*}
\prob{\theta}{f \in B} & = \int_{\mathcal{X}} f_B(x)
dp_{\theta}(x)  & \mbox{for\ } B\in
\mathcal{B} \mbox{\ \rv{20}}
\end{align*}
where 
\begin{itemize}
{\item $f_{\mathbb{R}}=\mathbf{1}$\ \rv{21}}
{\item $0\leq f_B \leq 1$\ \rv{22}}
{\item For countable disjoint $B_i$, $\sum_i f_{B_i} = f_{\bigcup_i B_i}$\
\rv{23}.}
\end{itemize}

As a remark, the result $f(x)$ is essentially a label. We could write the same
formula for functions with values in other measure spaces $(\mathcal{Y},\mathcal{B})$ than $\mathbb{R}$. Just let 
$\mathcal{B}$ be the $\sigma$-algebra on this space. In this way, we retrieve in
particular estimators in $\mathbb{R}^d$.

Another very important remark is that if we have access to two statistics
$f$ and $g$, we
have access to both \rv{24}. Indeed suppose that $f$ was taking its values in
$(\mathcal{Y},\mathcal{B})$ and $g$ in $(\mathcal{Z},\mathcal{C})$. Then take a new statistic with values in the
product space $(\mathcal{Y}\otimes\mathcal{Z},\mathcal{B}\otimes\mathcal{C})$, characterized by $h_{B\otimes C} =
f_B * g_C$ as real functions on $(\mathcal{X},\mathcal{A})$. We see that the three conditions are satisfied, and that the marginals
of $h$ are $f$ and $g$.

\subsubsection{From classical to quantum}
\label{transition}

The above description was already somewhat non-conventional, with the parallel
with quantum formalism in mind. In this subsection, we take one further step,
by setting classical probability as a special case of what will be our quantum
probability theory. 

To have something easy to understand, we start from a finite probability space
$(\mathcal{X},\mathcal{A})= \{1,\ldots,d\}$\ \rv{2}. We associate to it the Hilbert space of
complex valued functions on this space, that is $\mathcal{H} = \mathbb{C}^d$\ \rv{2}. We are here
endowed with a distinguished orthonormal basis $\{|e_i\rangle\}_{1\leq i\leq d}$
with $|e_i\rangle$ the function whose value is one on $i$ and zero elsewhere.

Notice by the way the notation $|\psi\rangle$: this is a physicist's notation for
vectors, elements of $\mathcal{H}$. They call this a ``ket''. The associated linear
form, that is, the adjoint of the vector, is called a ``bra'' and denoted
$\langle\psi|$. Thus $\langle\phi|\psi\rangle$ is the scalar product of
$|\phi\rangle$ and $|\psi\rangle$ (a ``bracket''). 

Now to the probability measure $p=(p_1,\ldots,p_d)$\ \rv{1} on $\{1,\ldots,d\}$, we
associate the matrix $\rho$\ \rv{1} diagonal in our special orthonormal basis
\rv{6}, with diagonal
entries $(p_1,\ldots,p_d)$.  As this is a diagonal matrix in an
orthonormal basis, with non-negative elements, this is a self-adjoint \rv{4}
non-negative \rv{5} matrix. Moreover, as $\sum_i p_i = 1$\ \rv{3}, it has trace
$1$\ \rv{3}. 

We see that the extremal points of our set are of matrices are the
orthogonal projectors on the lines spanned by our special eigenvectors, that is $|e_i\rangle\langle
e_i|$\ \rv{8}. They correspond to the Dirac measures on $i$. We may represent
any of these \emph{pure states}
by the eigenvector $|e_i\rangle$\ \rv{9}. We may also rewrite $\rho = \sum_i p_i
|e_i\rangle\langle e_i|$.

A statistical \emph{model}\ \rv{10} consists in a set of non-negative matrices
$\rho_{\theta}$ with trace $1$, on a Hilbert space $\mathcal{H}$, diagonal in the
$\{|e_i\rangle\}_i$ basis, indexed by a  parameter
$\theta$, for $\theta\in \Theta$\ \rv{11} the \emph{parameter space}. A
statistical problem consists in determining as precisely as possible, with a
meaning depending on the instance, a function of $\theta$.

As we have done for probability measures, we identify $f\in
L^{\infty}(\{1,\ldots,d\})$\ \rv{12,7} with the diagonal
matrix $O\in M(\mathbb{C}^d)$\ \rv{12,7}
whose diagonal elements are the $O_{i,i}= f(i)$. This is still the dual of
the set of matrices diagonal on our special basis. We view the action of $\rho$ by taking the trace of the
product with $\rho$. That is $p(f) = \xtr(\rho O)$\ \rv{15}. One can see that we
have only rewritten the classical formula  for the expectation.

Equivalently, measuring an observable $O$ yields as a result an eigenvalue of
$O$\ \rv{13}. The law of the result is given by:

\begin{align*}
\prob{\rho}{O \in B} & = \xtr(\rho P_{O,B}) & \mbox{for\ } B\in
\mathcal{B} \mbox{\ \rv{14}}
\end{align*}
where $P_{O,B}$ is the projection upon the space spanned by the eigenspaces of
$O$ corresponding to those eigenvalues $\lambda$ of $O$ such that $\lambda \in B$. In
other words, in our case, $O = \sum_i f(i) |e_i\rangle \langle e_i|$. Then
$P_{O,B}= \sum_{i| f(i)\in B} |e_i\rangle \langle e_i|$. This $P_{O,B}$ is
playing the role of ${\bf 1}_{f(x)\in B}$ in the classical setting. And we take
note that $\xtr(\rho P_{O,B}) = \sum_{i| f(i)\in B} p_i$, as we should obtain
from the classical formula.

We can encode in the same framework the general strategies for estimators, provided that
$\mathcal{X}_a$ is also finite\ \rv{16}. The auxiliary space is then identified to
$\mathcal{H}_a = \mathbb{C}^{d_a}$. We have
matrices $\rho_{\theta}\otimes\sigma$\ \rv{19}, with $\sigma$ independent of $\theta$.
We are allowed to use as observable $O$\ \rv{17} any matrix diagonal in the same basis as
these $\rho_{\theta}\otimes\sigma$.
The procedure equivalent to the partial integration on $\mathcal{X}_a$ is then taking
partial trace on $\mathcal{H}_a$ in
$\mathbb{P}_{\theta}[O\in  B] = \xtr((\rho_{\theta}\otimes\sigma)P_{O,B})$. And this yields
$\xtr(\rho_{\theta} M(B))$\ \rv{20} with
\begin{itemize}
{\item $M(\mathbb{R})=\mathbf{1}_{\mathcal{H}}$\ \rv{21}}
{\item $M(B)$ is non-negative and diagonal in the $\{|e_i\rangle\}$ basis
\rv{22}}
{\item For countable disjoint $B_i$, $\sum_i M(B_i) = M(\bigcup_i B_i)$\
\rv{23}.}
\end{itemize}

Here again, we see that if we have access to $O_1$ and $O_2$ characterized by
the families $M_1(B)$ and $M_2(C)$, we have access to both\ \rv{24}. Our new measurement
would be characterized by $N(B\otimes C) = M_1(B) M_2(C)$ as multiplication of
matrices. Notice that this set of matrices still satisfies the three above
conditions. Especially, the fact that they are still non-negative stems from
that they are diagonal in the same eigenbasis.

\smallskip

Going from classical to quantum now means throwing away our special eigenbasis
$\{|e_i\rangle\}$. The immediate consequence will be that we shall deal with
objects that do not commute. And of course, we did not restrain to
finite probability spaces in the classical case. Likewise,  we do
not restrain to finite-dimensional Hilbert spaces in the quantum case. We shall
therefore deal with operators rather than matrices. Keeping the
finite-dimensional example firmly in mind should be a guide to the intuition of
those less proficient in operator theory.

\subsubsection{Quantum}
\label{quantique}

A quantum system is described by a \emph{density operator} $\rho$\ \rv{1} over a
Hilbert space $\mathcal{H}$\ \rv{2}, that is:
\begin{dfntn}{\bf: Density operator}
\label{densop}

  A density operator, usually denoted by $\rho$, is a trace-class linear operator on a
(complex, separable) Hilbert space $ \mathcal{H} $ that satisfies:
\begin{itemize}
\item{$\rho$ is self-adjoint \rv{4}.}
\item{$\rho$ is non-negative (notice that this implies self-adjointness) \rv{5}.}
\item{$\xtr \,\rho =1$\ \rv{3}.}
\end{itemize}  
\end{dfntn}  

If $\mathcal{H}$ is finite-dimensional, those are just the (self-adjoint) non-negative
matrices with trace $1$.

We denote by $\mathcal{S}(\mathcal{H})$ the set of density operators on $\mathcal{H}$.

Density operators are a convex set, too. The extremal points are called ``pure
states''. They are the orthogonal projectors on $1$-dimensional spaces \rv{8}. Thus we can
represent them by a norm $1$ element of $\mathcal{H}$, denoted by $|\psi\rangle$ \rv{9}. The corresponding density matrix is  then $\rho= |\psi\rangle
\langle\psi|$. Notice that it would be more precise to speak of $|\psi\rangle$
as an element of the projective space $\mathcal{P}\mathcal{H}$, but we conform here to the
usage of physicists. Notice also that there are infinitely many pure states
even in the finite-dimensional case, unlike in the classical framework. Let us
finally signal that the decomposition of a mixed state on pure states is \emph{not}
unique. It is essentially unique if we further impose that the pure states of
the decomposition are all orthogonal, though.

A quantum statistical \emph{model}\ \rv{10} consists in a set of density operators
$\rho_{\theta}$ on a Hilbert space $\mathcal{H}$ indexed by a  parameter
$\theta$, for $\theta\in \Theta$\ \rv{11} the \emph{parameter space}. A
statistical problem consists in determining as precisely as possible, with a
meaning depending on the instance, a function of $\theta$.

Now the role of random variables is played by observables. Those are the
elements $O$\ \rv{12} of $\mathcal{B}_{sa}(\mathcal{H})$\ \rv{7}, the bounded self-adjoint operators upon
$\mathcal{H}$. If we are dealing with finite-dimensional $\mathcal{H}$, those are the
self-adjoint matrices. 

As a remark, the dual of $\mathcal{B}_{sa}(\mathcal{H})$ is the set of self-adjoint trace-class
operators, which $\rho$ is in. This duality is given by the formula of the expectation of
measuring $O$ on $\rho$, also called \emph{Born's rule}:
\begin{align}
\label{Born}
\mathbb{E}_{\rho}[O] = \xtr(\rho O) \ \ \ \mbox{\rv{15}}
\end{align} 

When measuring $O$, the result is an element of the spectrum of $O$\
\rv{13}, that
is in the finite-dimensional picture, an eigenvalue of $O$. The law of the
result when measuring $O$ on $\rho$ is:
\begin{align}
\label{mesproj}
\prob{\rho}{O \in B} & = \xtr(\rho P_{O,B}) & \mbox{for\ } B\in
\mathcal{B} \mbox{\ \rv{14}}
\end{align}
where $P_{O,B}$ is coming from the spectral measure of $O$. This is an object
associated to self-adjoint operators through the spectral theorem, whose main
property is that the expectation of the law above is given by the Born's rule
for any density operator $\rho$. We only give the derivation for finite-dimensional $\mathcal{H}$. Then, as $O$ is
self-adjoint, we can diagonalize it in an orthonormal basis, and write
$O=\sum_i \lambda_i |\psi_i\rangle \langle\psi_i|$. Then $P_{O,B}= \sum_{i|
\lambda_i\in B} |\psi_i\rangle \langle \psi_i|$. We see that in this case the law
of the measurement is coherent with the expectation given by Born's rule
(\ref{Born}).

Generally $\{P_{O,B}\}_B$ is a projector valued measure, the definition of
which we give below.  To each projector valued measure corresponds an
observable, and to each observable corresponds a projector valued measure. We
may then consider that this concept is also a definition of an observable. 
\begin{dfntn}{\bf: Projector valued measure\ \rv{12}}
\label{pvm}

  A projector operator valued measure $\{P(B)\}_{B\in \mathcal{B}}$ is a set of
operators on $\mathcal{H}$ such that:
\begin{itemize}
\item{$P(B)$ is an orthogonal projector.}
\item{$P(\mathbb{R}) = \bf{1}_{\mathcal{H}}$.}
\item{For disjoint countable $B_i$, $\sum_i P(B_i) = P(\bigcup_i B_i)$.}
\end{itemize}
\end{dfntn}
  
Notice that these are the axioms of a probability measure, except that we do not deal
with real numbers but with projection operators.

Combining this definition with the definition of a density operator, we can
check that formula 
(\ref{mesproj}) yields a true probability measure. Indeed, as both $\rho$ and
$P_{O,B}$ are non-negative, the probability of any event is non-negative. With the
countable additivity property of projector valued measure and linearity of product
and trace, we get the countable additivity of a probability measure. Finally, the
probability of the universe is $\xtr(\rho P_{O,\mathbb{R}}) = \xtr(\rho {\bf{1}_{\mathcal{H}}})
= 1$.

\emph{Remark:} - even for a pure state, the result of the measurement is random,
unless the pure state is an eigenvector of $O$.

 Now what is the most general estimation strategy, or measurement? The right analogy is that of
the auxiliary space. We measure observables $O$\ \rv{17} on a Hilbert space
$\mathcal{H}\otimes\mathcal{H}_a$\ \rv{18} under the density operator
$\rho_{\theta}\otimes\sigma$\ \rv{19}, with $\sigma$ independent of $\theta$. Now we may
take partial trace in (\ref{mesproj}) along $\mathcal{H}_a$, and we obtain equivalence
of this scheme with measuring a \emph{positive operator valued measure} (POVM).
\begin{dfntn}{\bf: Measurement (POVM) \rv{17}}
\label{measurement}

  A measurement M on a quantum system, taking values x in a measurable  space
$( \mathcal{X} , \mathcal{A} )$ is specified by a {\em positive operator valued
probability measure} or {\em POVM} for short, that is a collection of self-adjoint matrices $M(A): A \in \mathcal{A} $ such that:
\begin{itemize}
\item{$M( \mathcal{X} )= \mathbf{1}$, the identity matrix \rv{21}}
\item{Each $M(A)$ is non-negative \rv{22}}
\item{For disjoint countable $A_i$, $\sum_i M(A_i) = M(\bigcup A_i)$\ \rv{23}.}
\end{itemize}
The $M(A)$ are called the \emph{POVM elements}.
\end{dfntn}

The law of measuring $M$ on $\rho$ is given by
\begin{align}
\label{coll}
\prob{\rho}{O \in A} & = \xtr(\rho M(A)) & \mbox{for\ } A\in
\mathcal{A} \mbox{\ \rv{20}}.
\end{align}

  With the same reasoning as for projector valued measure (which are a special
case of these POVMs), this is a genuine probability measure.

A special case of POVM is that of a POVM dominated by $\sigma$-finite measure
$\nu$
on $(\mathcal{X},\mathcal{A})$, that is 
\begin{align}
\label{domination}
M(A) & = \int_{A} m(x)\dd\nu(x) \mbox{\ for all\ } A \in \mathcal{A}
\end{align}
where $m(x)$ is positive for all $x$ and $\int_{\mathcal{X}} m(x)\dd\nu(x) =
\mathbf{1}_{\mathcal{H}}$. The POVM associated to homodyne tomography is dominated by
the Lebesgue measure.

The very important difference with the classical world is that if we can have
access to $M_1$ or $M_2$, in general, we cannot have access to both
simultaneously\ \rv{24}. We cannot copy what we have done in the former paragraph,
since $M_1(A)M_2(B) + M_2(B)M_1(A)$ might not be non-negative if $M_1(A)$ and $M_2(B)$
do not commute. More generally, there is usually
no way to create a new POVM $N$ with values in $(\mathcal{X}\otimes\mathcal{Y},\mathcal{A}\otimes\mathcal{B})$ such
that the marginals are $M_1$ and $M_2$. Notably, two observables that do not
commute can never be measured simultaneously. As an example, consider that $M_1$ and
$M_2$ are two projector valued measures on $\mathbb{C}^2$, each with values in
$\{0,1\}$, corresponding to
observables diagonal in different bases $\{e_0,e_1\}$ and $\{f_0,f_1\}$. Then
$N(0,0)$ should be proportional both to $|e_0\rangle \langle e_0|$ and
$|f_0\rangle \langle f_0|$. So that it is null. Same remark for the other
$N(i,j)$. Thus $N(\{O,1\}^{\otimes 2}) = 0 \neq \mathbf{1}$. So that it is
null.

\smallskip

The truly quantum feature of quantum statistics lies in that we should decide
which measurement is to be carried out. Once we have chosen our measurement, we
are  left through (\ref{coll}) with a classical statistical experiment. This is
the case in this article.

As a last remark on the subject, we could have developed a slightly more
general formalism, based on $C^*$-algebras, that would have been parallel to Le
Cam formulation of statistics. In practical applications, the formalism above
is usually sufficient.

\subsection{Quantum homodyne tomography}
\label{QHD}

  The system we work with is the harmonic oscillator. Both in classical
or quantum mechanics, the harmonic oscillator is  a basic and pervading system. It
describes, notably, a particle on a  line, or a mode of the electromagnetic
field (that is monochromatic light), as in our case.

The state of a quantum harmonic oscillator is described by an operator on
$L^2(\mathbb{R})$ (this is the Hilbert space \rv{1}). There are two important
observables corresponding to the canonical coordinates of the particle. If
we know the expectation of measuring on a state $\rho$ any operator in the
algebra they generate, then we know $\rho$. Those observables
are  $ \mathbf{P} $,the magnetic field,    and $ \mathbf{Q}$, the electric
field. They satisfy the (canonical) commutation relations:

 \begin{equation*}
\begin{split}
[\mathbf{Q},\mathbf{P}] & = \mathbf{Q}\mathbf{P} - \mathbf{P}\mathbf{Q} \\ 
                        & =i \mathbf{1}. 
\end{split} 
\end{equation*}    

They are realized as:
\begin{eqnarray}
(\mathbf{Q} \psi_1)(x) & = & x \psi_1(x)  \notag \\
\label{expr}
(\mathbf{P} \psi_2)(x) & = & - i \frac{d\psi_2(x)}{dx}.  
\end{eqnarray}

As they do not commute, they cannot be measured simultaneously. However, any linear combination can theoretically be measured. These $ \mathbf{X}_\phi = \sin(\phi) \mathbf{Q} + \cos(\phi) \mathbf{P} $ are called {\em quadratures}.

  Using an experimental setup proposed in \cite{Vog},
  each of these quadratures could be experimentally measured on a laser beam \cite{Smithey}. The technique is called {\em quantum homodyne tomography}.

\begin{figure}[ht]
\centering
\input{QHDexp.pstex_t}
\caption{\label{QHDexp} Quantum Homodyne Tomography measurement set-up}
\end{figure}
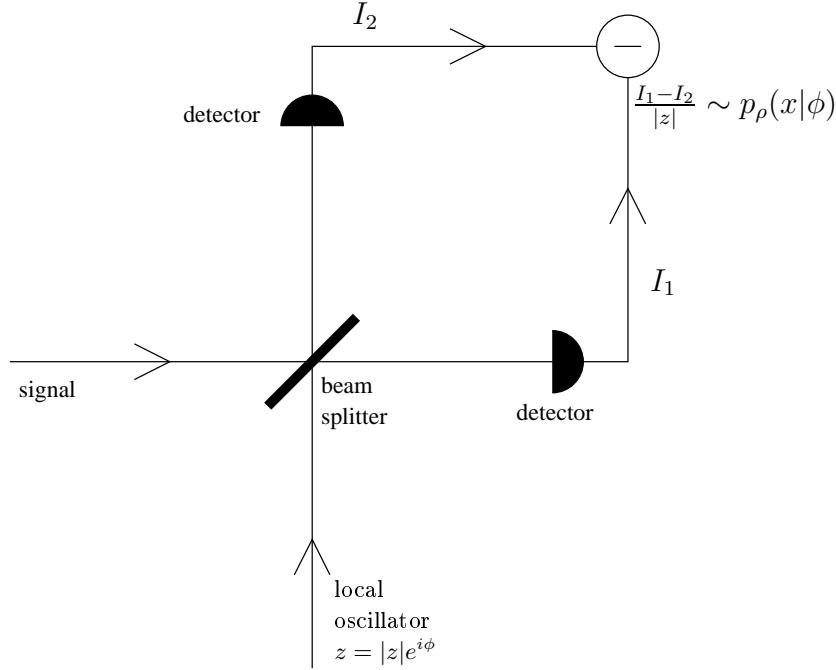

The optical set-up sketched in figure \ref{QHDexp} consists of an additional laser of high intensity $|z| \gg 1$ called the local oscillator,
a beam splitter through which the cavity pulse prepared in state $\rho$  is mixed with the
laser, and two photodetectors each measuring one of the two beams and producing
currents $I_{1,2}$ proportional to the number of photons. An electronic device produces the
result of the measurement by taking the difference of the two currents and rescaling it
by the intensity $|z|$.
    A simple quantum optics computation in \cite{Leonhardt}  shows that if the relative
phase between the laser and the cavity pulse is chosen to be $\phi$  then $(I_1 - I_2)/|z|$ has
density $p_{\rho}(x |\phi)$ corresponding to measuring $ \textbf{X}_{\phi}$ .

  Knowledge of $P_\rho(x|\phi)$, the law of the result of the measurement
$\mathbf{X}_\phi$ on $\rho$, for all $\phi$, is enough to reconstruct the state
$\rho$. As we have seen, the experimentalist may choose $\phi$ when measuring.
We assume that the measurement carried out on each of the $n$ systems in state $\rho$ is the following: first choose $\phi$ uniformly at random, then measure $ \mathbf{X}_\phi$. We get a random variable $ \mathbf{Y} =(\mathbf{X},\phi)$ with values in $ \mathbb{R} \times [0,\pi)  $ whose density with respect to the Lebesgue measure is $p_\rho(x,\phi) = \frac 1\pi p_\rho(x|\phi)$.

  Now we make explicit the links between $\rho$, $p_\rho(x,\phi)$ and the Wigner function $W_\rho$. First we write $\rho$ in a particular basis, physically very meaningful, the {\em Fock basis}, already given in Sec. \ref{Problème}: 
\begin{equation*}
\begin{split}
\psi_k(x) = H_k(x)e^{-x^2/\!2},
\end{split} 
\end{equation*} 
where $H_k$ is the $k$-th Hermite polynomial, normalized so that the $L^2$-norm of $\psi_k$ is $1$.
The projector on $\psi_k$ is the pure state with precisely $k$ photons. We  also denote this state by the ket $|k\rangle$.

The matrix entries of $p_\rho$ in this basis are $\rho_{j,k}=
\langle\psi_j,\rho \psi_k\rangle$. We can then derive from (\ref{Born}) and
(\ref{expr}) the formula we gave in Sec. \ref{Problème}:
\begin{eqnarray}
\bold{T}: \mathcal{S}( L^2(\mathbb{R})) & \longrightarrow     & L^1(\mathbb{R}\times [0,\pi]) \notag \\
\label{T}
            \rho   & \mapsto  &   \left(p_{\rho}:(x,\phi)\mapsto
\sum_{j,k=0}^{\infty}\rho_{j,k}\psi_j(x)\psi_k(x)e^{-i(j-k)\phi}\right).
\end{eqnarray}
 
  The mapping \textbf{T} associating $P_\rho$ to $\rho$ is invertible, so
we may hope to find $\rho$ from the independent identically distributed results $Y_1,Y_2,\dots,Y_n$ of the
measurements of the $n$ systems in state $\rho$. This implies notably that $p_\rho$ is another representation of the state.

  More explicitly,   there are {\em pattern functions} $f_{j,k}$  \cite{D'Ariano.0} against which to integrate $p_\rho$ to find any matrix entry of $\rho$ in the Fock basis, that is:
\begin{eqnarray*}
\rho_{j,k}=\int^\infty_{-\infty} dx \int_0^\pi \frac{d\phi}{\pi}p_\rho(x,\phi)
f_{j,k}(x)e^{i(j-k)\phi}.
\end{eqnarray*}

These $f_{j,k}$ are bounded real functions. That inverting the Radon transform
is an ill-posed problem can be seen in the behaviour of $f_{j,k}$ when $j$ and
$k$ go to infinity.  Several formulas were found for these functions \cite{D'Ariano.3}, among which:
\begin{equation}
\label{pattern}
 f_{j,k}(x)=\frac{d}{dx}(\chi_j(x)\phi_k(x))
\end{equation}
for $k\geq j$, where $\chi_j$ and $\phi_k$ are respectively the square-integrable and the unbounded solutions of the Schrödinger equation:
\begin{equation*} \left[ -\frac1{2} \frac{d^2}{dx^2} + \frac1{2} x^2 \right] \psi = \omega \psi,
\quad \omega \in \mathbb{R}. 
\end{equation*}
  
Another one, maybe more practical when it comes to theoretical calculations, or
when we add noise (see section \ref{noise}) is:
\begin{equation*}
  f_{j,k}(x,\phi)=\sqrt{\frac{j!}{k!}}\int_{-\infty}^{\infty}
|r|e^{-\frac{r^2}{2}+2 irx}r^{k-j}L_{j}^{k-j}(r^2) dr
\end{equation*}
where the $L_j^d$ are the Laguerre polynomials, that is the orthogonal
polynomials with respect to the measure $e^{-x}x^d$ on $\mathbb{R}^+$.

  Let's now have a look at the Wigner function. This is a real function of two
variables, with integral 1, but that may be negative in places. It can be
interpreted
as a generalized joint probability density of the electric and magnetic fields $q$ and
$p$. As both cannot be measured simultaneously, the negative patches are not nonsense. On the other hand, any projection on a line of the Wigner function
must be a true probability density, as it is the law of $ \mathbf{X}_\phi$,
which is an observable. In fact, the Wigner function may be seen as the
probability density on $\mathbb{R}^2$ resulting from (\ref{coll}) when measuring on $\rho$ a
``POVM'' whose elements are not non-negative, but whose marginals on each
line $\mathbb{R}$ are the $X_{\phi}$. 

As we have already said in the introduction, $p_\rho$ is the Radon transform of the Wigner function.  The Wigner function can be defined by its Fourier transform. 
This definition tells how to find the Wigner function $W$ of the state from its density matrix $\rho$: 
\begin{equation}
\label{rhoWigner}
\mathcal{F}_2 W(u,v) = \xtr(\rho e^{-iu\bold{Q}-iv\bold{P}}).
\end{equation}

On the other hand,  the generating function of $p_\rho(\cdot|\phi)$ is
\begin{eqnarray*}
\esp{e^{itX_\phi}} & = &  \xtr(\rho e^{it\bold{X}_\phi}).
\end{eqnarray*}
In other words,  $\mathcal{F}_2 W(t \cos\phi,t\sin\phi) = \mathcal{F}[p_\rho(\cdot,\phi)](t)$. These relations are known to imply that  $p_{\rho}  = \bold{R}(W)$ \cite{Deans} where \textbf{R} is the Radon transform. Explicitly:
\begin{eqnarray*}
p_{\rho}(x,\phi) & = & \int_{-\infty}^{\infty} W(x\cos\phi + y \sin \phi,x \sin
\phi - y\cos \phi)dy.
\end{eqnarray*} 
The Radon transform is illustrated by Fig. \ref{Radon}, given in Sec. \ref{Problème}.

  Finding the Wigner function from the data means then inverting the Radon transform, hence the name of tomography: that is the same mathematical problem as with the brain imagery technique called Positron Emission Tomography.

\subsection{Physical origin of the photocounter calibration problem}
\label{phy.photo}

  An experiment usually ends with a measurement. We need, however, an
apparatus to measure. And we first have to know what is the meaning of the
result the apparatus is giving us: it is not at all obvious a priori that if our
new thermometer says ``$31^\circ$ C'', the temperature cannot be ``$32^\circ$ C''. That is why
we must {\em calibrate} our measurement apparatus. In quantum mechanics, this
means associating with each result $i$ of our measurement the positive operator
$P(i)$, such that $P$ is the POVM (see definition \ref{measurement}) corresponding to our measurement. 

 In \cite{dariano-2004-93}, a general calibration procedure was intoduced. The procedure relies on comparing
with an already calibrated apparatus, using entangled states. Let
us describe this more precisely in the special case of the photocounter. 

A photocounter
is an apparatus that aims at counting the photons in a beam. The
ideal detector $D$ has therefore POVM elements given by $D(i)= |i\rangle\langle i|$
in the Fock basis. Recall we use  the physicists' notation, where
$|\cdot\rangle$ is a vector and $\langle\cdot|$ is the associated linear form. Moreover $|i\rangle$ is the  vector corresponding to the pure state
with $i$ photons, that is the function $\psi_i$ on $L^2( \mathbb{R})$, that  we
had defined in (\ref{hermite}).  

\begin{figure}
\centering
\input{photocounter3.pstex_t}
\caption{\label{photo} Experimental set-up to determine the POVM associated to
an unknown photocounter \textbf{P}. We use it to measure a known bipartite
state $|s\rangle$, jointly with a tomographer \textbf{T}. The photocounter gives a result
$i$ and the tomographer a result $x$. From these samples, we construct an estimator $\{\hat{P_i}\}$ of the self-adjoint operators associated to the results $\{i\}$ by the photocounter \textbf{P}.}
\end{figure}
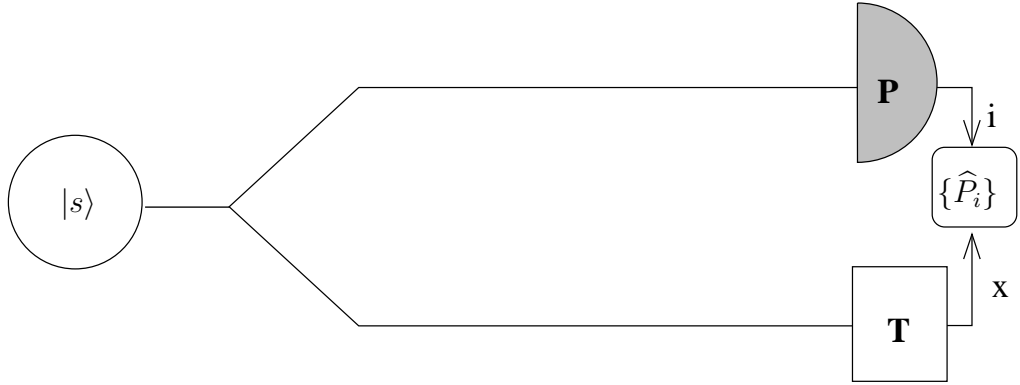

Models of the noise (non-unit efficiency and dark current)
leave the POVM diagonal in this basis. Thus, we are only interested in the
diagonal elements of $P_i$ in the Fock basis. To obtain those we send a twin beam state, one of the beams in the photocounter, the other in a homodyne tomographer. We get a result $i$ from the photo-counter, and $x$ from the tomographer (figure \ref{photo}; as we are only interested in the diagonal elements, we shall see that we do not need the phase $\phi$, as long as the experimentalist chooses it randomly). We then have to process these outcomes $(i,x)$ to find $P$.

   Mathematically, the twin beam is a system in a state
$|s\rangle=\sum_{k=0}^{\infty} b_k |k\rangle\otimes |k\rangle$. This notation
(where we may choose the $b_k$ non-negative) means that the
underlying Hilbert space is $L^2( \mathbb{R} )\otimes L^2( \mathbb{R})$, and
that $\rho$ is the pure state that projects on the line spanned by this vector.
Here again, $|k\rangle$ is the vector corresponding to the pure state
with $k$ photons. Finally $\sum_k b_k^2 =1$, so that the vector state $|s\rangle$ is normalized and the density
operator is $\rho = |s\rangle\langle s|$.

 Now, what is the law $p(i,x)$ of the samples we get? By (\ref{T}) we see that
the POVM associated to the tomographer is dominated by the Lebesgue measure on
$ \mathbb{R} \times [0,\pi)$, as in (\ref{domination}). That is $ {\langle
j|t_{x,\phi}|k\rangle}= \psi_j(x) \psi_k(x) e^{-i(j-k)\phi}$, where we have
denoted $t_{x,\phi}$ the self-adjoint operator associated to the result
$(x,\phi)$ for the POVM of the tomographer. If we forget about $\phi$ after
having chosen it randomly, we then get $\langle j|t_x|k\rangle = \psi_k(x)^2
\mathbf{1}_{j=k}$. We have now all the ingredients for calculating our law,
given the notation $\langle k|M_i|k\rangle = M_i^k$.

\begin{align}
  p(i,x) & = \xtr(\rho (P_i \otimes t_x)) \notag \\
         & = \langle s|(P_i \otimes t_x)|s\rangle \notag \\
         & = \sum_{k_1,k_2} b_{k_1}b_{k_2} (\langle k_1| \otimes\langle
k_1|)(P_i \otimes t_x)(|k_2\rangle \otimes |k_2\rangle) \notag \\
         & = \sum_{k_1,k_2} b_{k_1}b_{k_2} \langle k_1|P_i|k_2\rangle
\langle k_1|t_x|k_2\rangle \notag \\
         & = \sum_{k=0}^{\infty} b_k^2 P_i^k \psi_k(x)^2. \notag
\end{align} 

  (As a remark, the fourth line  shows that the use of the phase would be to retrieve the non-diagonal elements, in which we are not interested.)  

We have thus recovered \eqref{lawp}, and explained how we got the data with which we want to estimate the $M_i^m$.

\bibliographystyle{plain}  
\bibliography{bib_qlan}

\end{document}

%% file: QHDexp.pstex_t
\begin{picture}(0,0)%
\includegraphics{QHDexp.pstex}%
\end{picture}%
\setlength{\unitlength}{4144sp}%
\begingroup\makeatletter\ifx\SetFigFont\undefined%
\gdef\SetFigFont#1#2#3#4#5{%
  \reset@font\fontsize{#1}{#2pt}%
  \fontfamily{#3}\fontseries{#4}\fontshape{#5}%
  \selectfont}%
\fi\endgroup%
\begin{picture}(3905,4019)(1384,-4449)
\put(3338,-4011){\makebox(0,0)[lb]{\smash{{\SetFigFont{9}{10.8}{\rmdefault}{\mddefault}{\updefault}{\color[rgb]{0,0,0}local}%
}}}}
\put(3338,-4198){\makebox(0,0)[lb]{\smash{{\SetFigFont{9}{10.8}{\rmdefault}{\mddefault}{\updefault}{\color[rgb]{0,0,0}oscillator}%
}}}}
\put(3338,-4384){\makebox(0,0)[lb]{\smash{{\SetFigFont{9}{10.8}{\rmdefault}{\mddefault}{\updefault}{\color[rgb]{0,0,0}$z=|z|e^{i\phi}$}%
}}}}
\put(3444,-581){\makebox(0,0)[lb]{\smash{{\SetFigFont{12}{14.4}{\rmdefault}{\mddefault}{\updefault}{\color[rgb]{0,0,0}$I_2$}%
}}}}
\put(5226,-2176){\makebox(0,0)[lb]{\smash{{\SetFigFont{12}{14.4}{\rmdefault}{\mddefault}{\updefault}{\color[rgb]{0,0,0}$I_1$}%
}}}}
\put(5120,-1112){\makebox(0,0)[lb]{\smash{{\SetFigFont{12}{14.4}{\rmdefault}{\mddefault}{\updefault}{\color[rgb]{0,0,0}$\frac{I_1 - I_2}{|z|}\sim p_{\rho}(x|\phi)$}%
}}}}
\end{picture}%

%% file: photocounter3.pstex_t
\begin{picture}(0,0)%
\includegraphics{photocounter3.pstex}%
\end{picture}%
\setlength{\unitlength}{4144sp}%
\begingroup\makeatletter\ifx\SetFigFontNFSS\undefined%
\gdef\SetFigFontNFSS#1#2#3#4#5{%
  \reset@font\fontsize{#1}{#2pt}%
  \fontfamily{#3}\fontseries{#4}\fontshape{#5}%
  \selectfont}%
\fi\endgroup%
\begin{picture}(6053,2285)(4904,-3943)
\put(5223,-2918){\makebox(0,0)[lb]{\smash{{\SetFigFontNFSS{12}{14.4}{\rmdefault}{\mddefault}{\updefault}{\color[rgb]{0,0,0}$|s\rangle$}%
}}}}
\put(10468,-2858){\makebox(0,0)[lb]{\smash{{\SetFigFontNFSS{12}{14.4}{\rmdefault}{\mddefault}{\updefault}{\color[rgb]{0,0,0}$\{\widehat{P_i}\}$}%
}}}}
\end{picture}%